\documentclass{article}
\usepackage{times,mathptmx,amsmath,amsthm}

\newcommand{\bq}{\begin{equation}}
\newcommand{\eq}{\end{equation}}
\newcommand{\R}{{\bf R}}
\newcommand{\C}{{\bf C}}
\newcommand{\Z}{{\bf Z}}

\newcommand{\opname}[1]{\mathop{\fam0#1}}
\newcommand{\AD}{\opname{AD}}

\newcommand{\isumstar}{\opname{\sum\nolimits^*}\limits}
\newcommand{\sumstar}{\opname{\:\sum\nolimits^*}\limits}

\newcommand{\Tmin}{{T_{\rm min}}}
\newcommand{\Tmax}{{T_{\rm max}}}
\newcommand{\tTvert}{{\tilde{T}_{\rm vert}}}
\newcommand{\Podd}{{P^{\rm odd}}}
\newcommand{\Pnodd}{{P_n^{\rm odd}}}
\newcommand{\Peven}{{P^{\rm even}}}
\newcommand{\Pneven}{{P_n^{\rm even}}}

\newcommand{\cA}{{\cal A}}
\newcommand{\G}{{\rm G}}
\newcommand{\Grass}{{\rm Grass}}
\newcommand{\sgn}{{\rm sgn}}

\newcommand{\AzDi}{1}           
\newcommand{\RotMov}{2}         
\newcommand{\TiCo}{3}           
\newcommand{\StOr}{4}           
\newcommand{\HeFu}{5}           
\newcommand{\Conv}{6}           
\newcommand{\Two}{7}            
\newcommand{\Stacka}{8}         
\newcommand{\Stackb}{9}         
\newcommand{\Stackc}{10}        
\newcommand{\Stackd}{11}        
\newcommand{\Jag}{12}           
\newcommand{\EvenOdd}{13}       
\newcommand{\Rweight}{14}       
\newcommand{\Shuffle}{15}       
\newcommand{\Lemma}{16}         
\newcommand{\Lemmb}{17}         
\newcommand{\Ice}{18}           
\newcommand{\VeCo}{19}          
\newcommand{\Tilt}{20}          
\newcommand{\Trans}{21}         
                                                                                
\newcommand{\Ba}{{\bf 1}}       
\newcommand{\Be}{{\bf 2}}       
\newcommand{\CL}{{\bf 3}}       
\newcommand{\FW}{{\bf 4}}       
\newcommand{\Gr}{{\bf 5}}       
\newcommand{\Kaa}{{\bf 6}}      
\newcommand{\Kab}{{\bf 7}}      
\newcommand{\Li}{{\bf 8}}       
\newcommand{\Lo}{{\bf 9}}       
\newcommand{\MRRa}{{\bf 10}}    
\newcommand{\MRRb}{{\bf 11}}    
\newcommand{\Pe}{{\bf 12}}      
\newcommand{\PS}{{\bf 13}}      
\newcommand{\Ro}{{\bf 14}}      
\newcommand{\RR}{{\bf 15}}      
\renewcommand{\Sp}{{\bf 16}}      
\newcommand{\Sta}{{\bf 17}}     
\newcommand{\St}{{\bf 18}}      
\newcommand{\Th}{{\bf 19}}      
\newcommand{\Ya}{{\bf 20}}      

\begin{document}
 
\title{Alternating sign matrices and domino tilings}
\author{Noam Elkies\thanks{Supported by
                the Harvard Society of Fellows
                and by NSF grant DMS-87-18965.} \\
        Harvard College \\
    \\
        Greg Kuperberg\thanks{Supported by
                an NSF Graduate Research Fellowship.} \\
        University of California at Berkeley \\
    \\
        Michael Larsen\thanks{Supported by
                NSF grants DMS-8610730 and DMS-8807203.} \\
        University of Pennsylvania \\
    \\
        James Propp\thanks{Supported by
                an NSF Postdoctoral Research Fellowship.} \\
        Massachusetts Institute of Technology
    \\ }
 
\date{June 1991}

\maketitle
 
\begin{abstract}
We introduce a family of planar regions, called Aztec diamonds, and study the
ways in which these regions can be tiled by dominoes. Our main result is a
generating function that not only gives the number of domino tilings of the
Aztec diamond of order $n$ but also provides information about the orientation
of the dominoes (vertical versus horizontal) and the accessibility of one
tiling from another by means of local modifications. Several proofs of the
formula are given. The problem turns out to have connections with the
alternating sign matrices of Mills, Robbins, and Rumsey, as well as the square
ice model studied by Lieb.
\end{abstract}

\section{Introduction}

The {\bf Aztec diamond of order $n$} is the union of those lattice squares
$[a,a+1] \times [b,b+1] \subset \R^2$ ($a,b \in \Z$) that lie completely inside
the tilted square $\{(x,y): |x|+|y| \leq n+1\}$. (Figure~\AzDi \ shows the
Aztec diamond of order 3.) A {\bf domino} is a closed $1 \times 2$ or $2 \times
1$ rectangle in $\R^2$ with corners in $\Z^2$, and a {\bf tiling} of a region
$R$ by dominoes is a set of dominoes whose interiors are disjoint and whose
union is $R$. In this article we will show that the number of domino tilings of
the Aztec diamond of order $n$ is $2^{n(n+1)/2}$. We will furthermore obtain
more refined enumerative information regarding two natural statistics of a
tiling: the number of vertical tiles and the ``rank'' of the tiling (to be
defined shortly).

Fix a tiling $T$ of the Aztec diamond of order $n$. Every horizontal line $y=k$
divides the Aztec diamond into two regions of even area; it follows that the
number of dominoes that straddle the line must be even. Letting $k$ vary, we
see that the total number of vertical dominoes must be even; accordingly, we
define $v(T)$ as half the number of vertical tiles in $T$.

The most intuitively accessible definition of the rank-statistic $r(T)$ comes
by way of the notion of an ``elementary move'', which is an operation that
converts one domino-tiling of a region into another by removing two dominoes
that form a $2 \times 2$ block and putting them back rotated by 90 degrees (see
Figure~\RotMov). It will be shown that any domino-tiling of an Aztec diamond
can be reached from any other by a sequence of such moves; we may therefore
define the {\bf rank} of the tiling $T$ as the minimum number of moves required
to reach $T$ from the ``all-horizontals'' tiling (shown on the left side of
Figure~\RotMov). Thus the all-horizontals tiling itself has rank 0, while the
tiling shown  on the right side of Figure~\RotMov \ (viewed as a tiling of the
order 1 Aztec diamond) has rank 1.

Let 
$$\AD(n;x,q) = \sum_T x^{v(T)} q^{r(T)}$$
where $T$ ranges over all domino tilings of the order-$n$ Aztec diamond; this
is a polynomial in $x$ and $q$. The main result of this paper is:

\vspace{0.2in}

{\sc Theorem:}
$$\AD(n;x,q) = \prod_{k=0}^{n-1} (1+xq^{2k+1})^{n-k}.$$

\vspace{0.2in}

As important special cases, we have
\begin{eqnarray*}
\AD(n;q) & = & \prod_{k=0}^{n-1} (1+q^{2k+1})^{n-k}, \\
\AD(n;x) & = & (1+x)^{n(n+1)/2}, \ \mbox{and} \\
& & \\
\AD(n)   & = & 2^{n(n+1)/2},
\end{eqnarray*}
where we adopt the convention that an omitted variable is set equal to 1.

We will give four ways of understanding the formula for $\AD(n)$. The first
exploits the relationship between tilings of the Aztec diamond and the still
fairly mysterious ``alternating sign matrices'' introduced by Mills, Robbins,
and Rumsey in [\MRRa]. Our second proof yields the formula for $\AD(n)$ as a
special case of a theorem on monotone triangles (combinatorial objects closely
related to alternating sign matrices and introduced in [\MRRb]). The third
proof comes from the representation theory of the general linear group. The
last proof yields the more general formula for $\AD(n;x,q)$, and also leads to
a bijection between tilings of the order-$n$ diamond and bit-strings of length
$n(n+1)/2$. We conclude  by pointing out some connections between our results
and the ``square ice'' model studied in statistical mechanics.

\section{Height functions}

It is not at all clear from the definition of rank given in section 1 just how
one would calculate the rank of a specific tiling; for instance, it happens
that the all-verticals tiling of the order-$n$ Aztec diamond has rank
$n(n+1)(2n+1)/6$ and that every other tiling has strictly smaller rank, but it
is far from obvious how one would check this. Therefore, we will now give a
more technical definition of the rank, and prove that it coincides with the
definition given above. We use the vertex-marking scheme described in [\Th]; it
is a special case of the ``boundary-invariants'' approach to tiling problems
introduced in [\CL].

It will be conceptually helpful to extend a tiling $T$ of the Aztec diamond to
a tiling $T^+$ of the entire plane, by tiling the complement of the Aztec
diamond by horizontal dominoes in the manner shown in Figure~\TiCo \ for $n=3$.
Let $G$ be the graph with vertices $\{(a,b) \in \Z^2: |a|+|b| \leq n+1\}$, and
with an edge between $(a,b)$ and $(a',b')$ precisely when $|a-a'|+|b-b'|=1$.
Color the lattice squares of $\Z^2$ in black-white checkerboard fashion, so
that the line $\{(x,y): x+y = n+1\}$ that bounds the upper right border of the
Aztec diamond passes through only white squares. Call this the {\bf standard}
(or {\bf even}) {\bf coloring}. Orient each edge of $G$ so that a black square
lies to its left and a white square to its right; this gives the {\bf standard
orientation} of the graph $G$, with arrows circulating clockwise around white
squares and counterclockwise around black squares. (Figure~\StOr \ shows the
case $n=3$.) Write $u \rightarrow v$ if $uv$ is an edge of $G$ whose standard
orientation is from $u$ to $v$. Call $v=(a,b)$ a {\bf boundary vertex} of $G$
if $|a|+|b|=$ $n$ or $n+1$, and let the {\bf boundary cycle} be the closed
zigzag path $(-n-1,0), (-n,0), (-n,1), (-n+1,1), (-n+1,2), ..., (-1,n), (0,n),$
\linebreak
$(0,n+1), (0,n), (1,n), ...,
(n+1,0), ..., (0,-n-1), ..., (-n-1,0)$.
Call the vertex $v=(a,b)$ {\bf even} if it is the upper-left corner of a white
square (i.e., if $a+b+n+1$ is even), and {\bf odd} otherwise, so that in
particular the four corner vertices $(-n-1,0), (n+1,0), (0,-n-1), (0,n+1)$ are
even.

If one traverses the six edges that form the boundary of any domino, one will
follow three edges in the positive sense and three edges in the negative sense.
Also, every vertex $v$ of $G$ lies on the boundary of at least one domino in
$T^+$. Hence if for definiteness one assigns ``height'' $0$ to the leftmost
vertex $(-n-1,0)$ of $G$, there is for each tiling $T$ a unique way of
assigning integer-valued heights $H_T (v)$ to all the vertices $v$ of $G$,
subject to the defining constraint that if the edge $uv$ belongs to the
boundary of some tile in $T^+$ with $u \rightarrow v$, then $H_T(v) = H_T(u) +
1$. The resulting function $H_T (\cdot)$ is characterized by two properties:
\begin{description}
\item[\rm(i)]	$H(v)$ takes on the successive values
		$0,1,2,...,2n+1,2n+2,2n+1,...,0,$
		$...,2n+2,...,0$
		as $v$ travels along the boundary cycle of $G$;
\item[\rm(ii)]	if $u \rightarrow v$, then $H(v)$ is either
		$H(u)+1$ or $H(u)-3$.
\end{description}
The former is clear, since every edge of the boundary cycle is part of the
boundary of a tile of $T^+$. To see that (ii) holds, note that if the edge $uv$
belongs to $T^+$ (i.e.\ is part of the boundary of a tile of $T^+$), then $H(v)
= H(u) + 1$, whereas if $uv$ does not belong to $T^+$ then it bisects a domino
of $T^+$, in which case we see (by considering the other edges of that domino)
that $H(v) = H(u) - 3$.

In the other direction, notice that every height-function $H(\cdot)$ satisfying
(i) and (ii) arises from a tiling $T$, and that the operation $T \mapsto H_T$
is reversible: given a function $H$ satisfying (i) and (ii), we can place a
domino covering every edge $uv$ of $G$ with $|H(u)-H(v)|=3$, obtaining thereby
a tiling of the Aztec diamond, which will coincide with the original tiling $T$
in the event $H=H_T$. Thus there is a bijection between tilings of the Aztec
diamond and height functions $H(\cdot)$ on the graph $G$ that satisfy (i) and
(ii). For a geometric interpretation of $H(\cdot)$, see [\Th].

Figure~\HeFu \ shows the height-functions corresponding to two special tilings
of the Aztec diamond, namely (a) the all-horizontal tiling $\Tmin$ and (b) the
all-vertical tiling $\Tmax$. Since $H_T(v)$ is independent of $T$ modulo 4, we
are led to define the {\bf reduced height} $$h_T (v) = (H_T(v) -
H_{\Tmin}(v))/4;$$ parts (c) and (d) of Figure~\HeFu \ show the reduced
height-functions of $\Tmin$ and $\Tmax$, respectively. Lastly, we define the
rank-statistic $$r(T) = \sum_{v \in G} h_T(v).$$

It is easy to check that if one performs an elementary rotation on a 2-by-2
block centered at a vertex $v$ (a ``$v$-move'' for short), the effect is to
leave $h_T (v')$ alone for all $v' \neq v$ and to either increase or decrease
$h_T (v)$ by 1; we call the move {\bf raising} or {\bf lowering} respectively.

We may now verify that $r(T)$ (as defined by the preceding equation) is equal
to the number of elementary moves required to get from $T$ to $\Tmin$. Since
$r(\Tmin) = 0$, and since an elementary move merely changes the reduced height
of a single vertex by $\pm 1$, at least $r(T)$ moves are required to get from
$T$ to $\Tmin$. It remains to check that for every tiling $T$ there is a
sequence of moves leading from $T$ to $\Tmin$ in which only $r(T)$ moves are
made. To find such a sequence, let $T_0 = T$ and iterate the following
operation for $i=0,1,2,...\,$: Select a vertex $v_i$ at which $h_{T_i}(\cdot)$
achieves its maximum value. If $h_{T_i}(v_i) = 0$, then $T_i = \Tmin$ and we
are done. Otherwise, we have $h_{T_i}(v_i) > 0$, so that $H_{T_i}(v_i) >
H_{\Tmin}(v_i)$, with $v_i$ not on the boundary of $G$ (since $h_{T_i}$
vanishes on the boundary). Reducing $H_{T_i}(v_i)$ by 4 preserves the legality
of the height-labelling, and corresponds to performing a $v_i$-move on $T_i$,
yielding a new tiling $T_{i+1}$ with $r(T_{i+1}) = r(T_i)-1$. By repeating this
process, we continue to reduce the rank-statistic by 1 until the procedure
terminates at $T_{r(T)} = T_{\min}$.

Thus we have shown that every tiling of the Aztec diamond may be reached from
every other by means of moves of the sort described. This incidentally
furnishes another proof that the number of dominoes of each orientation
(horizontal or vertical) must be even, since this is clearly true of $\Tmin$
and since every move annihilates two horizontal dominoes and creates two
vertical ones, or vice versa.

The partial ordering on the set of tilings of an Aztec diamond given by
height-functions has a pleasant interpretation in terms of a two-person game.
Let $T,T'$ be tilings of the Aztec diamond of order $n$. We give player A the
tiling $T$ and player B the tiling $T'$. On each round, A makes a rotation
move, and B has the choice of either making the identical move (assuming it is
available to her) or passing. Here, to make an ``identical move'' means to find
an identically-situated 2-by-2 block in the identical orientation and give it a
90-degree twist. If, after a certain number of complete rounds (i.e. moves by A
and counter-moves by B) A has solved her puzzle (that is, reduced the tiling to
the all-horizontals tiling) while B has not, then A is deemed the winner;
otherwise, B wins. Put $T' \preceq T$ if and only if B has a winning strategy
in this game. It is easily checked  (without even considering any facts about
tilings) that the relation $\preceq$ is reflexive, asymmetric, and transitive.
In fact, $T' \preceq T$ if and only if $h_T (v) \leq h_{T'} (v)$ for all $v \in
G$. Moreover, the ideal strategy for either player is to make only lowering
moves -- though in the case $T' \preceq T$, it turns out that B can win by
copying A whenever possible, regardless of whether such moves are lowering or
raising.

\section{Alternating sign matrices}

An {\bf alternating sign matrix} is a square matrix ($n$-by-$n$, say) all of
whose entries are $1$, $-1$, and $0$, such that every row-sum and column-sum is
1, and such that the non-zero entries in each row and column alternate in sign;
for instance
$$
\left( \begin{array}{rrrr}
	0  &\ 1  &\ 0  &\ 0  \\
	1  & -1  &\ 1  &\ 0  \\
	0  &\ 0  &\ 0  &\ 1  \\
	0  &\ 1  &\ 0  &\ 0 \end{array} \right)
$$
is a typical 4-by-4 alternating sign matrix. (For an overview of what is
currently known about such matrices, see [\Ro].) Let $\cA_n$ denote the set of
$n$-by-$n$ alternating sign matrices.

If $A$ is an $n$-by-$n$ alternating sign matrix with entries $a_{ij}$ ($1 \leq
i,j \leq n$), we may define
$$a_{ij}^* = i \ + \ j \ - 
\ 2 \left( \sum_{i'=1}^i \sum_{j'=1}^j a_{i'j'} \right)$$
for $0 \leq i,j \leq n$. We call the $(n+1)$-by-$(n+1)$ matrix $A^*$ the {\bf
skewed summation} of $A$. (It is a variant of the ``corner-sum matrix'' of
[\RR].) The matrices $A^*$ that arise in this way are precisely those such that
$a_{i0}^* = a_{0i}^* = i$ and $a_{in}^* = a_{ni}^* = n-i$ for $0 \leq i \leq
n$, and such that adjacent entries of $A'$ in any row or column differ by 1.
Note that $a_{ij} =  \frac12(a_{i-1,j}^* + a_{i,j-1}^* - a_{i-1,j-1}^* -
a_{i,j}^*)$, so that an alternating sign matrix can be recovered from its
skewed summation. Thus, the alternating sign matrix $A$ defined above has
$$
A^* =
\left( \begin{array}{rrrrr}
	0 & 1 & 2 & 3 & 4 \\
	1 & 2 & 1 & 2 & 3 \\
	2 & 1 & 2 & 1 & 2 \\
	3 & 2 & 3 & 2 & 1 \\
	4 & 3 & 2 & 1 & 0 \end{array} \right)
$$
as its skewed summation.

Our goal is to show that the domino tilings of the Aztec diamond of order $n$
are in 1-to-1 correspondence with pairs $(A,B)$ where $A \in \cA_n$, $B \in
\cA_{n+1}$, and $A,B$ jointly satisfy a certain ``compatibility'' relation. We
will do this via the height-functions defined in the previous section.

Given a tiling $T$ of the order-$n$ Aztec diamond, we construct matrices $A'$
and $B'$ that record $H_T (v)$ for $v$ odd and even, respectively (where
$v=(x,y) \in G$ is even or odd according to the parity of $x+y+n+1$). We let
$$ a'_{ij} = H_T (-n+i+j,-i+j) $$
for $0 \leq i,j \leq n$ and
$$ b'_{ij} = H_T (-n-1+i+j,-i+j) $$
for $0 \leq i,j \leq n+1$; thus, the tiling of Figure~\Conv \ gives the
matrices
$$
A' = 
\left( \begin{array}{rrrrr}
	1 & 3 & 5 & 7 & 9 \\
	3 & 5 & 7 & 5 & 7 \\
	5 & 7 & 5 & 7 & 5 \\
	7 & 5 & 3 & 5 & 3 \\
	9 & 7 & 5 & 3 & 1 \end{array} \right)
\ \mbox{and} \ 
B' =
\left( \begin{array}{rrrrrr}
	 0 & 2 & 4 & 6 & 8 & 10 \\
	 2 & 4 & 6 & 4 & 6 &  8 \\
	 4 & 6 & 8 & 6 & 4 &  6 \\
	 6 & 4 & 6 & 4 & 6 &  4 \\
	 8 & 6 & 4 & 2 & 4 &  2 \\
	10 & 8 & 6 & 4 & 2 &  0 \end{array} \right) .
$$
Note that the matrix-elements on the boundary of $A'$ and $B'$ are independent
of the particular tiling $T$. Also note that in both matrices, consecutive
elements in any row or column differ by exactly 2. Therefore, under suitable
normalization, $A'$ and $B'$ can be seen as skewed summations of alternating
sign matrices $A$ and $B$. Specifically, by putting $a_{ij}^* = (a'_{ij} -
1)/2$ and $b_{ij}^* = b'_{ij}/2$, we arrive at matrices $A^*,B^*$ which, under
the inverse of the skewed summation operation, yield the matrices $A,B$ that we
desire:
$$
A^* =
\left( \begin{array}{rrrrr}
	0 & 1 & 2 & 3 & 4 \\
	1 & 2 & 3 & 2 & 3 \\
	2 & 3 & 2 & 3 & 2 \\
	3 & 2 & 1 & 2 & 1 \\
	4 & 3 & 2 & 1 & 0 \end{array} \right)
\ \mbox{and} \ 
B^* =
\left( \begin{array}{rrrrrr}
	0 & 1 & 2 & 3 & 4 & 5 \\
	1 & 2 & 3 & 2 & 3 & 4 \\
	2 & 3 & 4 & 3 & 2 & 3 \\
	3 & 2 & 3 & 2 & 3 & 2 \\
	4 & 3 & 2 & 1 & 2 & 1 \\
	5 & 4 & 3 & 2 & 1 & 0 \end{array} \right) ;
$$
$$
A = \left( \begin{array}{rrrr}
	0  &  0  &  1  &  0 \\
	0  &  1  & -1  &  1 \\
	1  &  0  &  0  &  0 \\
	0  &  0  &  1  &  0 \end{array} \right)
\ \mbox{and} \ 
B = \left( \begin{array}{rrrrr}
	0  &  0  &  1  &  0  &  0 \\
	0  &  0  &  0  &  1  &  0 \\
	1  &  0  &  0  & -1  &  1 \\
	0  &  1  &  0  &  0  &  0 \\
	0  &  0  &  0  &  1  &  0 \end{array} \right) .
$$
Conversely, $A$ and $B$ determine $A'$ and $B'$, which determine $H_T$, which
determines $T$.

There is an easy way of reading off $A$ and $B$ from the domino-tiling $T$,
without using height-functions. First, note that the even vertices in the
interior of the Aztec diamond of order $n$ are arranged in the form of a tilted
$n$-by-$n$ square. Also note that each such vertex is incident with 2, 3, or 4
dominoes belonging to the tiling $T$; if we mark each such site with a $1$,
$0$, or $-1$ (respectively), we get the entries of $A$, where the upper-left
corner of each matrix corresponds to positions near the left corner of the
diamond. Similarly, the odd vertices of the Aztec diamond (including those on
the boundary) form a tilted $(n+1)$-by-$(n+1)$ square. If we mark each such
site with a $-1$, $0$, or $1$ according to whether it is incident with 2, 3, or
4 dominoes of the extended tiling $T^+$, we get the entries of $B$. (We omit
the proof that this construction agrees with the one we gave earlier, since it
is only the first one that we actually need.)

The legality constraint (ii) from the previous section tells us that for $1
\leq i,j \leq n$, the internal entries $b'_{ij}$ of the matrix $B'$ must be
equal to
$$
\begin{array}{llll}
\mbox{either} & a'_{i-1,j-1} - 1     
& \ \mbox{or} & a'_{i-1,j-1} + 3, \\
\mbox{either} & a'_{i-1,j}   - 3   
& \ \mbox{or} & a'_{i-1,j}   + 1, \\
\mbox{either} & a'_{i,j-1}   - 3   
& \ \mbox{or} & a'_{i,j-1}   + 1, \ \mbox{and} \\
\mbox{either} & a'_{i,j}     - 1 
& \ \mbox{or} & a'_{i,j}     + 3 .
\end{array}
$$
Thus, in all but one of the six possible cases for the submatrix
$$\left(
\begin{array}{cc}
a'_{i-1,j-1}   & a'_{i-1,j} \\
a'_{i,j-1} & a'_{i,j}
\end{array}
\right)$$
shown in Table~\Two, the value of $b'_{ij}$ is uniquely determined; only in the
case
$$\left(
\begin{array}{cc}
2k-1 & 2k+1 \\
2k+1 & 2k-1
\end{array}
\right)$$
arising from $a_{ij} = 1$ does $b'_{ij}$ have two possible values, namely
$2k-2$ and $2k+2$.

\setcounter{table}{\Two}
\addtocounter{table}{-1}
\begin{table}
\begin{center}
\begin{tabular}{l|ccc|c}
$\left( \begin{array}{ll}
        a'_{i-1,j-1} & a'_{i-1,j} \\
        a'_{i,j-1}   & a'_{i,j} \end{array} \right)$
& & $a_{ij}$ & & $b'_{ij}$ \\
& & & & \\
\hline
& & & & \\
$\left( \begin{array}{ll} 2k-1 & 2k+1 \\ 2k+1 & 2k+3 \end{array} \right)$
& & 0 & & $2k+2$ \\
& & & & \\
$\left( \begin{array}{ll} 2k+3 & 2k+1 \\ 2k+1 & 2k-1 \end{array} \right)$
& & 0 & & $2k+2$ \\
& & & & \\
$\left( \begin{array}{ll} 2k+1 & 2k-1 \\ 2k+3 & 2k+1 \end{array} \right)$
& & 0 & & $2k$ \\
& & & & \\
$\left( \begin{array}{ll} 2k+1 & 2k+3 \\ 2k-1 & 2k+1 \end{array} \right)$
& & 0 & & $2k$ \\
& & & & \\
$\left( \begin{array}{ll} 2k-1 & 2k+1 \\ 2k+1 & 2k-1 \end{array} \right)$
& &  1 & & $2k-2 \ \mbox{or} \ 2k+2$ \\
& & & & \\
$\left( \begin{array}{ll} 2k+1 & 2k-1 \\ 2k-1 & 2k+1 \end{array} \right)$
& & -1 & & $2k$
\end{tabular}
\end{center}
\caption{Two-by-two submatrices}
\end{table}

It now follows that if we hold $A$ fixed, the number of $(n+1)$-by-$(n+1)$
alternating sign matrices $B$ such that the pair $(A,B)$ yields a legal height
function is equal to $2^{N_+(A)}$, where $N_+(A)$ is the number of $+1$'s in
the $n$-by-$n$ alternating sign matrix $A$. That is:
\begin{equation}
\label{AD1}
\AD(n) = \sum_{A \in \cA_n} 2^{N_+(A)} \ .
\end{equation}

Switching the roles of $A$ and $B$, we may by a similar argument prove
\begin{equation}
\label{AD2}
\AD(n) = \sum_{B \in \cA_{n+1}} 2^{N_-(B)} \ ,
\end{equation}
where $N_-(\cdot)$ gives the number of $-1$'s in an alternating sign matrix.
Replacing $n$ by $n-1$ and $B$ by $A$ in (\ref{AD2}), we get
\begin{equation}
\label{AD3}
\AD(n-1) = \sum_{A \in \cA_n} 2^{N_-(A)} \ .
\end{equation}
On the other hand, $N_+(A) = N_-(A)+n$ for all $A \in \cA_n$, so (\ref{AD1})
tells us that
\begin{equation}
\label{AD4}
\AD(n) = 2^n \sum_{A \in \cA_n} 2^{N_-(A)} \ .
\end{equation}
Combining (\ref{AD3}) and (\ref{AD4}), we derive the recurrence relation
$$\AD(n) = 2^n \AD(n-1),$$ which suffices to prove our formula for $\AD(n)$.
(Mills, Robbins and Rumsey [\MRRa] prove 
$$\sum_{A \in \cA_n} 2^{N_-(A)} = 2^{n(n-1)/2}$$
as a corollary to their Theorem 2.)

In the remainder of this section, we discuss tilings and alternating sign
matrices from the point of view of lattice theory. Specifically, we show that
the tilings of an order-$n$ Aztec diamond correspond to the lower ideals (or
``down-sets'') of a partially ordered set $P_n$, while the $n$-by-$n$
alternating sign matrices correspond to the lower ideals of a partially ordered
set $Q_n$, such that $P_n$ consists of a copy of $Q_n$ interleaved with a copy
of $Q_{n+1}$. (For terminology associated with partially ordered sets, see
[\Sta].)

We start by observing that the set of legal height functions $H$ on the
order-$n$ Aztec diamond is a poset in the obvious component-wise way, with $H_1
\geq H_2$ if $H_1(v) \geq H_2(v)$ for all $v \in G$. Moreover, the consistency
conditions (i) and (ii) are such that if $H_1$ and $H_2$ are legal
height-functions, then so are $H_1 \vee H_2$ and $H_1 \wedge H_2$, defined by
$(H_1 \vee H_2)(v) = \max(H_1(v),H_2(v))$ and $(H_1 \wedge H_2)(v) =
\min(H_1(v),H_2(v))$; thus our partially ordered set is actually a distributive
lattice.

$\cA_n$, the set of $n$-by-$n$ alternating sign matrices, also has a lattice
structure. Given $A_1,A_2 \in \cA_n$, we form their skewed summations
$A_1^*,A_2^*$, and declare $A_1 \geq A_2$ if every entry of $A_1^*$ is greater
than or equal to the corresponding entry of $A_2^*$. This partial ordering on
alternating sign matrices is intimately connected with the partial ordering on
tilings: if $T_1,T_2$ are tilings, then $T_1 \geq T_2$ if and only if $A_1 \geq
A_2$  and $B_1 \geq B_2$, where $(A_i,B_i)$ is the pair of alternating sign
matrices corresponding to the tiling $T_i$ ($i = 1,2$).

For each vertex $v=(x,y)$ of the graph $G$ associated with the order-$n$ Aztec
diamond (with $x,y \in \Z$, $|x|+|y| \leq n+1$), let $m=H_\Tmin(v)$ and
$M=H_\Tmax(v)$, so that $m,m+4,...,M$ are the possible values of $H_T(v)$, and
introduce points  $(x,y,m), (x,y,m+4), ..., (x,y,M-4) \in \Z^3$ lying above the
vertex $v=(x,y)$. (Note that if $(x,y)$ is on the boundary of $G$, $m=M$, so
the set of points above $(x,y)$ is empty.) Let $P$ denote the set of all such
points as $v$ ranges over the vertex-set of $G$. We make $P$ a directed graph
by putting an edge from $(x,y,z) \in P$ to $(x',y',z') \in P$ provided $z=z'+1$
and $|x-x'|+|y-y'|=1$; we then make $P$ a partially ordered set by putting
$(x,y,z) \geq (x',y',z')$ if there is a sequence of arrows leading from
$(x,y,z)$ to $(x',y',z')$.

To each height-function $H$ we may assign a subset $I_H \subseteq P$, with
$I_H=\{(x,y,z) \in P: z < H(x,y)\}$. This operation is easily seen to be a
bijection between the legal height-functions $H$ and the lower ideals of the
partially ordered set $P$. Indeed, the natural lattice structure on the set of
height-functions $H$ (with $H_1 \leq H_2$ precisely if $H_1(v) \leq H_2(v)$ for
all $v \in G$) makes it isomorphic to the lattice $J(P)$ of lower ideals of
$P$, and the rank $r(T)$ of a tiling $T$ (as defined above) equals the rank of
$H_T$ in the lattice, which in turn equals the cardinality of $I_{H_T}$.

Note that for all $(x,y,z) \in P$, $x+y+z \equiv n+1$ mod 2. The poset $P$
decomposes naturally into two complementary subsets $\Peven$ and $\Podd$, where
a point $(x,y,z) \in P$ belongs to $\Peven$ if $z$ is even and $\Podd$ if $z$
is odd. The vertices of $\Peven$ form a regular tetrahedral array of side $n$,
resting on a side (as opposed to a face); that is, it consists of a $1$-by-$n$
array of nodes, above which lies a $2$-by-$(n-1)$ array of nodes, above which
lies a $3$-by-$(n-2)$ array of nodes, and so on, up to the $n$-by-$1$ array of
nodes at the top. The partial ordering of $P$ restricted to $\Peven$ makes
$\Peven$ a poset in its own right, with $(x,y,z)$ covering $(x',y',z')$ when
$z=z'+2$ and $|x-x'|=|y-y'|=1$. Similarly, the vertices of $\Podd$ form a
tetrahedral array of side $n-1$; each vertex of $\Podd$ lies at the center of a
small tetrahedron with vertices in $\Peven$. $\Podd$, like $\Peven$, is a poset
in itself, with $(x,y,z)$ covering $(x',y',z')$ when $z=z'+2$ and
$|x-x'|=|y-y'|=1$.

Our correspondence between height-functions $H_T$ and pairs $(A,B)$ of
alternating sign matrices tells us that $\cA_{n}$, as a lattice, is isomorphic
to $J(\Pnodd)$, while $\cA_{n+1}$ is isomorphic to $J(\Pneven)$. Indeed, under
this isomorphism, $A \in \cA_n$ and $B \in \cA_{n+1}$ are compatible if and
only if the union of the down-sets of $\Peven$ and $\Podd$ corresponding to $A$
and $B$ is a down-set of $P= \Peven \cup \Podd$. (This coincides with the
notion of compatibility given in [\RR].) If we let $Q_n$ denote the tetrahedral
poset $\Pnodd$ (so that $\Pneven$ is isomorphic to $Q_{n+1}$), then we see that
$P_n$ indeed consists of a copy of $Q_n$ interleaved with a copy of $Q_{n+1}$.

As an aid to visualizing the poset $P$ and its lower ideals, we may use stacks
of marked 2-by-2-by-2 cubes resting on a special multi-level tray. The bottom
face of each cube is marked by a line joining midpoints of two opposite edges,
and the top face is marked by another such line, skew to the mark on the bottom
face. (See Figure~\Stacka.) These marks constrain the ways in which we allow
ourselves to stack the cubes. To enforce these constraints, whittle away the
edges of the cube on the top and bottom faces that are parallel to the marks on
those faces, and replace each mark by a protrusion, as in Figure~\Stackb; the
rule is that a protrusion on the bottom face of a cube must fit into the space
between two whittled-down edges (or between a whittled-down edge and empty
space). The only exception to this rule is at the bottom of the stack, where
the protrusions must fit into special furrows in the tray. Figure~\Stackc \
shows the tray in the case $n=4$; it consists of four levels, three of which
float in mid-air. On the bottom level, the outermost two of the three gently
sloping parallel lines running from left to right should be taken as
protrusions, and the one in between should be taken as a furrow. Similarly, in
the higher levels of the tray, the outermost lines are protrusions and the
innermost two are furrows. We require that the cubes that rest on the table
must occupy only the $n$ obvious discrete positions; no intermediate positions
are permitted. Also, a cube cannot be placed unless its base is fully
supported, by the tray, a tray and a cube, or two cubes.

In stacking the cubes, one quickly sees that in a certain sense one has little
freedom in how to proceed; any stack one can build will be a subset of the
stack shown in Figure~\Stackd \ in the case $n=4$. Indeed, if one partially
orders the cubes in Figure~\Stackd \ by the transitive closure of the relation
``is resting on'', then the poset that results is the poset $P$ defined
earlier, and the admissible stacks correspond to lower ideals of $P$ in the
obvious way. Moreover, the markings visible to an observer looking down on the
stack yield a picture of the domino tiling that corresponds to that stack.

\section{Monotone triangles}

Let $A^*$ be the skewed summation of an $n$-by-$n$ alternating sign matrix $A$.
Notice that the $i$th row ($0 \leq i \leq n$) begins with an $i$ and ends with
an $n-i$, so that reading from left to right, we must see $i$ descents and
$n-i$ ascents; that is, there are exactly $i$ values of $j$ in $\{1,2,...,n\}$
satisfying $a'_{i,j} = a'_{i,j-1} - 1$, and the remaining $n-i$ values of $j$
satisfy $a'_{i,j} = a'_{i,j-1} + 1$. Form a triangular array whose $i$th row
($1 \leq i \leq n$) consists of those values of $j$ for which $a'_{i,j} =
a'_{i,j-1} - 1$. E.g., for
$$
A = \left( \begin{array}{rrrr}
	0 &  1 & 0 & 0 \\
	1 & -1 & 1 & 0 \\
	0 &  0 & 0 & 1 \\
	0 &  1 & 0 & 0
    \end{array} \right) \ \mbox{and} \ 
A' = \left( \begin{array}{ccccc}
	0 & 1 & 2 & 3 & 4 \\
	1 & 2 & 1 & 2 & 3 \\
	2 & 1 & 2 & 1 & 2 \\
	3 & 2 & 3 & 2 & 1 \\
	4 & 3 & 2 & 1 & 0   
     \end{array} \right)$$
we get the triangle
$$\begin{array}{ccccccc}
   &   &   &   &   &   &   \\
   &   &   & 2 &   &   &   \\
   &   &   &   &   &   &   \\
   &   & 1 &   & 3 &   &   \\
   &   &   &   &   &   &   \\
   & 1 &   & 3 &   & 4 &   \\
   &   &   &   &   &   &   \\
 1 &   & 2 &   & 3 &   & 4 \\
   &   &   &   &   &   &   \\
\end{array}$$
Note that $j$ occurs in the $i$th row of the monotone triangle exactly if the
sum of the first $i$ entries in column $j$ of the alternating sign matrix is 1.

A {\bf monotone triangle} of size $n$ is a triangular array of natural numbers
with strict increase from left to right along its $n$ rows and with non-strict
increase from left to right along its diagonals, as in the array above. If the
bottom row of a monotone triangle is $1 \ 2 \ ... \ n$, we call the array a
{\bf complete monotone triangle}. It is not hard to show that the preceding
construction gives a bijection between  the $n$-by-$n$ alternating sign
matrices and the complete monotone triangles of size $n$. Moreover, the $+1$'s
in the alternating sign matrix correspond to entries in some row of the
triangle that do not occur in the preceding row.

It follows from the foregoing that $\AD(n)$ is the sum, over all complete
monotone triangles of size $n$, of 2 to the power of the number of entries in
the monotone triangle that do not occur in the preceding row. Since a monotone
triangle of size $n$ has exactly $n(n+1)/2$ entries, we may divide both sides
of the equation $\AD(n) = 2^{n(n+1)/2}$ by $2^{n(n+1)/2}$ and paraphrase it as
the claim that the sum, over all complete monotone triangles of size $n$, of
$\frac12$ to the power of the number of entries in the monotone triangle that
{\it do} occur in the preceding row, is precisely 1.

Define the weight of a monotone triangle (of any size) as $\frac12$ to the
power of the number of entries that appear in the preceding row, and let
$W(a_1, a_2, ..., a_k)$ be the sum of the weights of the monotone triangles of
size $k$ with bottom row $a_1 \ a_2 \ ... \ a_k$. (For now we may assume $a_1 <
a_2 < ... < a_k$, although we will relax this restriction shortly.) Our goal is
to prove that $W(1,2,...,n) = 1$ for all $n$.

To this end, observe that we have the recurrence relation
\begin{equation}
\label{recrel}
W(a_1,a_2,...,a_n) =
\sumstar_{b_1=a_1}^{a_2}
\sumstar_{b_2=a_2}^{a_3}
...
\sumstar_{b_{n-1}=a_{n-1}}^{a_n}
W(b_1,b_2,...,b_{n-1})
\end{equation}
for all $n$, where $\isumstar$ is the modified summation operator
$$\sumstar_{i=r}^s f(i) =
\frac12 f(r) + f(r+1) + f(r+2) + ... + f(s-1) + \frac12 f(s) \ ;$$
for the number of factors of $\frac12$ that contribute to the coefficient of
$W(b_1,b_2,...,$ $b_{n-1})$ is exactly the number of $b_i$'s that also occur
among the $a_i$'s. Observe that the operator $\isumstar$ resembles definite
integration in that
\begin{equation}
\label{defint}
\sumstar_{i=r}^s f
+
\sumstar_{i=s}^t f
=
\sumstar_{i=r}^t f
\end{equation}
for $r < s < t$. Indeed, if we extend $\sumstar$ by defining 
$$\sumstar_{i=r}^r f = 0$$
for all $r$ and $$\sumstar_{i=r}^s f = - \sumstar_{i=s}^r f$$ for $r > s$, then
$(\ref{defint})$ holds for all integers $r,s,t$. Hence, starting from the
base-relation $W(a_1) = 1$, equation $(\ref{recrel})$ can be applied
iteratively to define $W(\cdot)$ as a function of $a_1,a_2,...,a_n$, regardless
of whether $a_1 < a_2 < ... < a_n$ or not.

Notice that if $f(x) = x^m$, then 
$$\sumstar_{r=s}^t f(r)$$
is a polynomial in $s$ and $t$ of degree $m+1$,
of the form
$$\frac{t^{m+1}-s^{m+1}}{m+1} \ + \ \mbox{terms of lower order.}$$
More generally, if $f(x,y,z,...)$ is a polynomial in $x,y,z,...$ with a
highest-order term $c x^{\!m} y^{m'} z^{m''} ...$, then $$\sumstar_{r=s}^t
f(r,y,z,...)$$ is a polynomial in $s,t,y,z,...$ of degree $\deg f + 1$,  with
highest-order terms
$$  \frac{c}{m+1} t^{m+1} y^{m'} z^{m''} \dots \ \ \  \mbox{and}
\ \ \ - \frac{c}{m+1} s^{m+1} y^{m'} z^{m''} \dots \ .$$

We will now use $(\ref{recrel})$ and $(\ref{defint})$
to prove the general formula
\begin{equation}
\label{genform}
W(a_1,a_2,...,a_n) = \prod_{1 \leq i < j \leq n} 
\frac{a_j - a_i}{j-i} \ .
\end{equation}
(This immediately yields $W(1,2,...,n) = 1$, which as we have seen implies
$\AD(n) = 2^{n(n+1)/2}$.) Formula (\ref{genform}) is equivalent to Theorem 2 in
[\MRRb], but we offer our own proof.

The proof is by induction. When $n=1$, we have $W(a_1) = 1$, so that
$(\ref{genform})$ is satisfied. Suppose now that we have
$$W(b_1,b_2,...,b_{n-1}) = \prod_{1 \leq i < j \leq n-1} 
\frac{b_j - b_i}{j-i}$$
for all $b_1, b_2, ..., b_{n-1}$. Since $W(b_1, b_2, ..., b_{n-1})$ is a
polynomial of degree 
\linebreak
$(n-1)(n-2)/2$ with a highest-order term
$$\frac{b_2^1 \, b_3^2 \, \cdots \, b_{n-1}^{n-2}}
{1! \, 2! \, \cdots \, (n-2)!} \ ,$$
the recurrence relation $(\ref{recrel})$
and the observations made in the preceding paragraph
imply that $W(a_1,a_2,...,a_n)$ is
a polynomial of degree
$$(n-1)(n-2)/2 + (n-1) = n(n-1)/2$$
with a highest-order term
$$\frac{a_2^1 \, a_3^2 \, \cdots \, a_{n}^{n-1}}
{1! \, 2! \, \cdots \, n!} \ .$$

To complete the proof,
we need only show that
$W(a_1, a_2, ..., a_n)$
is skew-symmetric in its arguments;
for this implies
that it is divisible by
\linebreak
$(a_2 - a_1) (a_3 - a_1) \cdots (a_n - a_{n-1})$,
a polynomial of the same degree
(namely 
\linebreak
$n(n-1)/2$) as itself,
and a comparison of the coefficients of leading terms
yields $(\ref{genform})$.

It suffices to show that
interchanging any two consecutive arguments of $W$
changes the sign of the result.
For convenience,
we illustrate with $n=4$:
\begin{eqnarray*}
W(a_2, a_1, a_3, a_4) & = &
\sumstar_{b_1=a_2}^{a_1}
\sumstar_{b_2=a_1}^{a_3}
\sumstar_{b_3=a_3}^{a_4}
W(b_1,b_2,b_3) \\
& = &
\left(-\sumstar_{b_1=a_1}^{a_2} \right)
\left(\sumstar_{b_2=a_1}^{a_2} + \sumstar_{b_2=a_2}^{a_3} \right)
\left(\sumstar_{b_3=a_3}^{a_4} \right)
W(b_1,b_2,b_3) \ .
\end{eqnarray*}
The skew-symmetry of $W(b_1,b_2,b_3)$
kills off one of the two terms:
\begin{eqnarray*}
& & \sumstar_{b_1=a_1}^{a_2}
    \sumstar_{b_2=a_1}^{a_2}
    \sumstar_{b_3=a_3}^{a_4}
    W(b_1,b_2,b_3) \\
& = & \sumstar_{b_2=a_1}^{a_2}
    \sumstar_{b_1=a_1}^{a_2}
    \sumstar_{b_3=a_3}^{a_4}
    W(b_2,b_1,b_3) \\
& & \mbox{(by re-labelling)} \\
& = & \sumstar_{b_1=a_1}^{a_2}
    \sumstar_{b_2=a_1}^{a_2}
    \sumstar_{b_3=a_3}^{a_4}
    W(b_2,b_1,b_3) \\
& & \mbox{(by commutativity)} \\
& = & - \sumstar_{b_1=a_1}^{a_2}
    \sumstar_{b_2=a_1}^{a_2}
    \sumstar_{b_3=a_3}^{a_4}
    W(b_1,b_2,b_3) \\
& & \mbox{(by skew-symmetry)} ,
\end{eqnarray*}
implying that the term vanishes.
Hence
\begin{eqnarray*}
W(a_2,a_1,a_3,a_4)
& = & - \sumstar_{b_1=a_1}^{a_2}
    \sumstar_{b_2=a_2}^{a_3}
    \sumstar_{b_3=a_3}^{a_4}
    W(b_1,b_2,b_3) \\
& = & -W(a_1,a_2,a_3,a_4) \ ,
\end{eqnarray*}
as claimed.
Similarly,
$W(a_1,a_2,a_4,a_3) = - W(a_1,a_2,a_3,a_4)$.
A slightly more complicated calculation,
involving a sum of four terms of which three vanish,
gives $W(a_1,a_3,a_2,a_4) = -W(a_1,a_2,a_3,a_4)$.
The argument for the skew-symmetry of
$W(a_1,a_2,...,a_n)$
is much the same for $n$ in general,
although the notation is more complex;
we omit the details.

Having shown that $W$ is skew-symmetric
in its arguments,
we have completed the proof of $(\ref{genform})$,
which yields the formula for $\AD(n)$
as a consequence.

Some remarks are in order.
First, it is noteworthy that
$$\prod_{1 \leq i < j \leq n}
\frac{a_j - a_i}{j-i}$$
is an integer
provided $a_1,...,a_n$ are;
this can be proved in 
a messy but straightforward manner
by showing that every prime $p$
must divide the numerator
at least as many times as
it divides the denominator.
Alternatively, one can show
that this product is equal to the
determinant of the $n$-by-$n$ matrix
whose $i,j$th entry is
the integer
$$a_i \choose {j-1}$$
(see [\PS] and [\Sp]).

Second, formula $(\ref{genform})$
has a continuous analogue:
If we take $V(x) = 1$
for all real $x$
and inductively define
$$V(x_1,x_2,...,x_n) =
\int_{x_1}^{x_2}
\int_{x_2}^{x_3}
\cdots
\int_{x_{n-1}}^{x_n}
V(y_1,y_2,...,y_{n-1}) \:
dy_{n-1} \:
\cdots \:
dy_2 \:
dy_1 \: ,$$
then essentially the same argument shows that
$$V(x_1,x_2,...,x_n) = \prod_{1 \leq i < j \leq n} \frac{x_j - x_i}{j-i} \ .$$
This has the following probabilistic interpretation:
Given $n$ real numbers $x_1 < x_2 < ... < x_n$,
let $X_{i,i} = x_i$ for $1 \leq i \leq n$,
and for all $1 \leq i < j \leq n$
let $X_{i,j}$ be a number chosen
uniformly at random in
the interval $[x_i,x_j]$.
Then the probability that $X_{i,j} \leq X_{i+1,j}$
and $X_{i,j} \leq X_{i,j+1}$ for all suitable $i,j$ is
$$\prod_{1 \leq i < j \leq n} \frac{1}{j-i}
= \frac{1}{1!\,2!\, \cdots \,(n-1)!} \ .$$
We do not know a more direct proof of this fact
than the one outlined here.

Third, the usual (unstarred) summation operator
does not satisfy a relation like $(\ref{defint})$,
so the method used here
will not suffice to count unweighted monotone triangles.
(Mills, Robbins, and Rumsey
offer abundant evidence that the number of complete monotone triangles
of size $n$ is
$$\prod_{k=0}^{n-1} \frac{(3k+1)!}{(n+k)!} \ ,$$
but no proof has yet been found.)
However,
the operators
$$\sum_{i=r}^s {\,\!}^L = \sum_{i=r}^{s-1}$$
and
$$\sum_{i=r}^s {\,\!}^R = \sum_{i=r+1}^{s}$$
{\it do} satisfy an analogue of $(\ref{defint})$,
and one can exploit this to give
streamlined proofs of some formulas
in the theory of plane partitions;
details will appear elsewhere.

Fourth,
we should note that the function $W(a_1,a_2,...,a_m)$
has significance for tilings of the Aztec diamond of order $n$,
even outside the case with $m=n$ and $a_i = i$ for $1 \leq i \leq m$.
Suppose $m \leq n$ and $a_m \leq n$,
and let $\Pi$ be the path in the graph $G$
that starts at $(-m,n-m)$
whose $2j-1$st and $2j$th steps head south and east respectively
if $j \in \{a_1,...,a_n\}$
and otherwise head east and south respectively,
for $1 \leq j \leq n$,
ending at the vertex $(n-m,-m)$;
Figure \Jag(a) shows $\Pi$
when $n=4$, $m=2$, $(a_1,a_2) = (2,3)$.
It is not hard to show that
the number of domino tilings of the
portion of the Aztec diamond
that lies above $\Pi$
is $2^{m(m+1)/2} W(a_1,a_2,...,a_m)$.

Fifth (and last),
we should note that
the role played by the matrix $A$
at the beginning of the section
(in expressing $\AD(n)$ in terms of
a weighted sum over complete monotone triangles
of size $n$)
could have been played just as well
by the matrix $B$,
giving rise to an alternative formula
expressing $\AD(n)$ as
the sum,
over all complete monotone triangles of size $n+1$,
of 2 to the power of
the number of entries above the bottom row
that do not occur in the succeeding row.
But, dividing by $2^{n(n+1)/2}$,
we reduce the claim $\AD(n) = 2^{n(n+1)/2}$
to the same claim as before
(the sum of the weights of
all the fractionally weighted complete monotone triangles
of any given size
is equal to 1).
This gives a second significance of $W(\cdots)$
for tilings of the Aztec diamond.
Specifically,
suppose $m \leq n+1$ and $a_m \leq n$,
and let $\Pi$ be the path in the graph $G$
that starts at $(-m,n+1-m)$
whose $2j-1$st and $2j$th steps
head east and south respectively
if $j \in \{a_1,...,a_n\}$
and otherwise head south and east respectively,
for $1 \leq j \leq n+1$,
ending at the vertex $(n+1-m,-m)$;
Figure \Jag(b) shows $\Pi$ when
$n=4$, $m=2$, $(a_1,a_2) = (2,3)$.
It can be shown that
the number of domino tilings
of the portion of the Aztec diamond
that lies above $\Pi$
is $2^{m(m-1)/2} W(a_1,a_2,...,a_m)$.

\section{Grassmann algebras}

The resemblance between the formula for $W$ and the Weyl dimension formula for
representations of $GL(n)$ is not coincidental. In fact, the identity
$$\sum_A 2^{N(A)} = 2^{n(n-1)/2}$$
can be proved by pure representation theory. The idea is to relate the rules
for consecutive rows in Gelfand triangles to the decomposition of
$GL(n)$-representations as $GL(n-1) \times GL(1)$ representations.

Let $V$ be a finite-dimensional vector space, $\Lambda^i (V)$ the $i$th
exterior power of $V$, and
$$\Grass(V) = \bigoplus_{i=0}^n \Lambda^i (V), $$
the Grassmann algebra generated by $V$.
It is elementary that if $V$ and $W$ are finite-dimensional vector spaces,
$$\Lambda^2 (V \oplus W) = \Lambda^2 (V) \oplus \Lambda^2 (W) \oplus
	(V \otimes W),$$
$$\Grass(V \oplus W) = \Grass(V) \otimes \Grass(W).$$
Writing $\G(V)$ for $\Grass(\Lambda^2(V))$, we get
$\dim(\G(V)) = 2^{\dim(V) \choose 2}$ and
$$\G(V \oplus W) = \G(V) \otimes \G(W) \otimes \Grass(V \otimes W).$$

We now recall the Cartan-Weyl theory of weights of irreducible representations
of Lie groups, in the case of $GL(n)$ (due to Schur); for more details, see
[\Gr]. If $V = \C^n$, then $GL(V) = GL(n,\C)$ contains the group $T$ of
diagonal matrices $diag(x_1,...,x_n)$. The analytic homomorphisms $T
\rightarrow \C^*$ are precisely the Laurent monomials $x_1^{a_1} \cdots
x_n^{a_n}$, $a_i \in \Z$. If $\rho$ is a finite-dimensional (analytic)
representation of $GL(n,\C)$, its restriction to $T$ is a direct sum of
1-dimensional analytic representations (called {\it weights\/}), and the
restriction of the trace of $\rho$ to $T$ is a Laurent polynomial in the $x_i$;
we represent a weight of $\rho$ by the sequence of exponents occurring in the
corresponding Laurent monomial in $tr(\rho|_T)$. For instance, the trace
function of the identity representation is the sum of the diagonal elements,
$x_1+\cdots+x_n$, so the weights are the basis vectors $(0,...,0,1,0,...,0)$.
The operations of linear algebra can be translated into operations on trace
polynomials. Thus, the trace of a direct sum of representations is the sum of
the traces, the trace of a tensor product is the product of traces, the trace
of the $k$th exterior power  is the $k$th elementary symmetric function of the
constituent monomials, and so on. The irreducible representations $\rho$ of
$GL(n,\C)$ are indexed by {\it dominant weights} $\lambda = (\lambda_1, ...,
\lambda_n)$, where $\lambda_i \in \Z$ and $\lambda_1 \leq ... \leq \lambda_n$;
among all weights occurring in $\rho|_T$ satisfying this inequality, $\lambda$
has the greatest norm. For instance, the dominant weight of the identity
representation is $(0,0,...,0,1)$.

We set $a_i = \lambda_i + i$, so the finite-dimensional irreducible
representations of $GL(n,\C)$ are indexed by $(a_1,...,a_n)$ where $a_i \in \Z$
and $a_1 < ... < a_n$. The Weyl character formula for $GL(n)$ says that the
trace function for $\rho(a_1,...,a_n)$ is
$$
\frac
{\sum_{\sigma \in S_n} \sgn(\sigma) \, 
	{x_1}^{a_{\sigma(1)}} \cdots {x_n}^{a_{\sigma(n)}}}
{\prod_{1 \leq i \leq n} x_i \ 
 \prod_{1 \leq j < i \leq n} (x_i - x_j)} .
$$
The numerator of this expression can be written
$$
\det
\left( \begin{array}{cccc}
	x_1^{a_1} & x_2^{a_1} & \dots & x_n^{a_1} \\
	x_1^{a_2} & x_2^{a_2} & \dots & x_n^{a_2} \\
	\vdots & \vdots & \ddots & \vdots \\
	x_1^{a_n} & x_2^{a_n} & \dots & x_n^{a_n}
       \end{array} \right) \ .
$$
Subtracting $x_1^{a_{i+1}-a_i}$ times row $i$ from row $i+1$,
for $i = n-1,n-2,...,1$, we obtain
$$x_1^{a_1} \sum_{\tau \in S_{n-1}} \sgn(\tau)
\prod_{i=1}^{n-1} \left( 
 x_{i+1}^{a_{\tau(i)+1}} - x_{i+1}^{a_{\tau(i)}} x_1^{a_{\tau(i)+1}-a_{\tau(i)}}
\right) . $$
The trace functions for the representations $V$, $\Lambda^2(V)$, and $\G(V)$
(that is, for the action of $GL(n)$ on these spaces induced by the action of
$GL(n)$ on $V$) are given by
$\sum_{1 \leq i \leq n} x_i$,
$\sum_{1 \leq j < i \leq n} x_i x_j$, and
$\sum_{1 \leq j < i \leq n} (1 + x_i x_j)$, respectively.
The trace function for $\rho(a_1,...,a_n) \otimes \G(V)$ is therefore
$$
x_1^{a_1} \sum_{\tau \in S_{n-1}} \sgn(\tau)
\prod_{i=1}^{n-1} \left( 
 x_{i+1}^{a_{\tau(i)+1}} - x_{i+1}^{a_{\tau(i)}} x_1^{a_{\tau(i)+1}-a_{\tau(i)}}
\right) 
\prod_{i=1}^n \frac{1}{x_i} 
\prod_{1 \leq j < i \leq n} \frac{1+x_i x_j}{x_i - x_j};
$$
this is equal to
$$x_1^{a_1-1} \sum_{\tau \in S_{n-1}} \sgn(\tau)
\prod_{i=1}^{n-1} (\Sigma_1 (i,\tau) + \Sigma_2 (i,\tau))
\prod_{i=2}^{n} \frac{1}{x_i} 
\prod_{2 \leq j < i \leq n} \frac{1+x_i x_j}{x_i - x_j},$$
where
\begin{eqnarray*}
\Sigma_1 (i,\tau) & = & \sum_{k=a_{\tau(i)}}^{a_{\tau(i)+1}-1}
	x_{i+1}^k x_1^{a_{\tau(i)+1}-1-k} \ , \\
\Sigma_2 (i,\tau) & = & \sum_{\ell=a_{\tau(i)}+1}^{a_{\tau(i)+1}}
	x_{i+1}^\ell x_1^{a_{\tau(i)+1}-1-\ell} \ .
\end{eqnarray*}
Viewing $GL(n-1,\C)$ as the subgroup of $GL(n,\C)$
consisting of all matrices of the form
$\left( \begin{array}{cc} 1 & 0 \\ 0 & M \end{array} \right)$,
we can restrict $\rho(a_1,...,a_n) \otimes \G(\C^n)$ to $GL(n-1)$.
At the level of traces on the diagonal, this amounts to setting $x_1=1$,
to obtain
\small
$$
\left( 2 \sum_{b_1 = a_1}^{a_2} {\!\! ^*} \right)
\cdots
\left( 2 \sum_{b_{n-1} = a_{n-1}}^{a_n} {\!\!\!\!\!\! ^* \ } \right)
\sum_{\tau \in S_{n-1}} \sgn(\tau)
x_2^{b_{\tau(1)}} \cdots x_n^{b_{\tau(n-1)}}
\prod_{i=2}^{n} \frac{1}{x_i} 
\prod_{2 \leq j < i \leq n} \frac{1+x_i x_j}{x_i - x_n}.
$$
\normalsize
(As in the preceding section, the notation $\sumstar$ indicates a sum where the
endpoints are counted with multiplicity $\frac{1}{2}$.) This is visibly the sum
of the traces of the $GL(n-1)$-representations $$\rho(b_1,...,b_{n-1}) \otimes
\G(\C^{n-1}),$$ counted with appropriate multiplicities. In fact, since two
representations of $GL(n-1)$ are the same if and only if their trace
polynomials coincide, this gives a formula for the restriction of the
representation $\rho(a_1,...,a_n) \otimes \G(\C^n)$ to $GL(n-1)$. Iterating
this process, we see that $W(1,2,...,n)$ (as defined in the previous section)
is the value obtained by substituting $x_1 = x_2 = ... = x_n = 1$ in the trace
function of $\G(\C^n)$ viewed as a $GL(n)$-representation, or in other words,
the trace function of $\G(\C^n)$ on $GL(0)=1$, which is simply the dimension of
$\G(\C^n)$.

\section{Domino shuffling}

The even (or standard) coloring of the Aztec diamond, as defined earlier, is
the black-white checkerboard coloring in which the interior squares along the
northeast border are black. In this section, it will be convenient to also
consider the other checkerboard coloring, which we call {\bf odd}. We will
continue to call a vertex of a checkerboard-colored region {\bf even} if it is
the upper-left corner of a white square and {\bf odd} otherwise, only now this
notion depends on the checkerboard coloring chosen as well as on the
coordinates of the vertex.

In general, a union of squares in a bi-colored checkerboard will be called {\it
even} if the leftmost square in its top row is white, and {\it odd} if that
square is black. Thus, the left half of Figure~\EvenOdd shows an even Aztec
diamond, an even 2-by-2 block, and two even dominoes (along with an even
vertex), while the right half of Figure~\EvenOdd shows odd objects of the same
kind.

Given a tiling of a colored region by dominoes, we may remove all the odd
blocks to obtain an {\bf odd-deficient tiling.} In general, an odd-deficient
domino tiling of a region in the plane is a partial tiling that has no odd
blocks and that can be extended to a complete tiling of that region by adding
only odd blocks. An odd-deficient tiling of the Aztec diamond of order $n$ with
its even coloring is uniquely determined by the heights of its even vertices,
as recorded in the matrix $B$ of section 3; thus, these odd-deficient tilings
are in one-one correspondence with alternating sign matrices of order $n+1$.

Given a partial tiling $\tilde{T}$ of the plane, let $U_{\tilde{T}}$ be the
union of the dominoes belonging to $\tilde{T}$. Observe that if $\tilde{T}$ is
odd-deficient, then the boundary of $U_{\tilde{T}}$ has corners only at odd
vertices.

The functions $v(T)$ and $r(T)$ defined earlier can be expressed in the form
$$v(T) = \sum_{d \in T} v(d)$$
and
$$r(T) = \sum_{d \in T} r(d)$$
for suitable functions $v(\cdot)$ and $r(\cdot)$ on the set of dominoes, which
we now define. If the domino $d$ is horizontal, let $v(d) = r(d) = 0$; if $d$
is vertical, let $v(d) = \frac12$ and let $r(d)$ be assigned according to the
location of the center of $d$ following the pattern set down in Figure~\Rweight
\ for the case $n=3$. (More formally, we may declare that if $d$ is the
vertical domino with upper-left corner at $(i,j)$, then $r(d) = (-1)^{i+j+n}
(i+n+1)$.) Clearly $v(T)$ is the sum of $v(d)$ over all dominoes $d \in T$. As
for $r(T)$, note that
$$r(\Tmin) = 0 = \sum_{d \in \Tmin} r(d) \ ;$$
also note that a move that increases $h(T)$ by $1$ either creates two vertical
dominoes $d_1,d_2$ satisfying $r(d_1)+r(d_2)=1$ or annihilates two vertical
dominoes $d_1,d_2$ satisfying $r(d_1)+r(d_2)=-1$. Thus by induction $r(T) =
\sum_{d \in T} r(d)$ for all tilings $T$.

We therefore have
$$\AD(n;x,q) = \sum_T \prod_{d \in T} x^{v(d)} q^{r(d)}.$$
We now prove
$$\AD(n;x,q) =
\prod_{k=0}^{n-1} (1+xq^{2k+1})^{n-k}$$
using a process called {\bf domino-shuffling}, which is a certain involution on
the set of odd-deficient tilings of an infinite checkerboard. If $d$ is domino
on a colored region, we define $S(d)$, the {\bf shuffle} of $d$, as the domino
obtained by moving $d$  one unit to the left or up if it is even  and one unit
to the right or down if it is odd.   (See Figure~\Shuffle.) Graphically, one
can put an arrow joining the two non-corner vertices on the boundary of $d$,
pointing from the even vertex to the odd vertex; this indicates the direction
in which $d$ will shuffle.

Clearly $S$ is an involution on the set of dominoes on an infinite
checkerboard. Two dominoes form an odd block if and only if each is the shuffle
of the other; if $d$ and $S(d)$ are horizontal, then $r(d) + r(S(d)) = 0$,
while if $d$ and $S(d)$ are vertical, then $r(d) + r(S(d)) = -1$.

Given a partial tiling $\tilde{T}$ we define $S(\tilde{T})$, the {\bf shuffle}
of $\tilde{T}$, to be the collection of all $S(d)$ with $d \in \tilde{T}$.

\vspace{0.2in}

{\sc Lemma:} Domino shuffling is an involution on the odd-deficient tilings of
an infinite checkerboard.

\vspace{0.2in}

Proof: Let $\tilde{T}$ be an odd-deficient tiling of the plane, with $T$ an
extension to a true tiling of the plane. We first show that $S(\tilde{T})$ is a
partial tiling, that is, that no two dominoes of $S(\tilde{T})$ overlap. Assume
otherwise, and suppose that a white square $s$ is covered by two dominoes in
$S(\tilde{T})$. That is, $S(\tilde{T})$ contains two of the four dominoes
$a,b,c,d$ shown in Figure~\Lemma \  (with arrows indicating the directions in
which they shuffle). There are six cases to be considered and ruled out.

$a,b \in S(\tilde{T})$: $\tilde{T}$ must contain the dominoes $S^{-1}(a) =
S(a)$ and $S^{-1}(b) = S(b)$. But $S(a)$ and $S(b)$ overlap (see
Figure~\Lemmb(a)).

$c,d \in S(\tilde{T})$: Same reasoning.

$a,c \in S(\tilde{T})$: $S(a), \, S(c) \in \tilde{T}$ (see Figure~\Lemmb(b)).
The full tiling $T$ must cover $s$ but cannot include $b$ or $d$ (since $T$
already includes $S(a)$ and $S(c)$ which conflict with those two dominoes);
hence $T$ must include $a$ or $c$. But in the former case, $a \in T$ forms an
odd block with $S(a) \in T$, so that $S(a) \not \in \tilde{T}$ after all; and
the case $c \in T$ leads to a similar contradiction.

$b,d \in S(\tilde{T})$: Same reasoning.

$a,d \in S(\tilde{T})$: Same reasoning as in the preceding two cases, though
the geometry is somewhat different (see Figure~\Lemmb(c)).

$b,c \in S(\tilde{T})$: Same reasoning.

Hence a white square cannot be covered by two dominoes of $S(\tilde{T})$. The
proof for black squares is similar. Therefore, $S(\tilde{T})$ is a partial
tiling of the checkerboard.

We must also show that $S(\tilde{T})$ is odd-deficient. $S(\tilde{T})$ cannot
contain any odd blocks, because the inverse shuffle (which is the same as the
shuffle) of an odd block is an odd block. It remains to show that the boundary
of $U_{S(\tilde{T})}$ has corners only at odd vertices. Let $v$ be an even
vertex. It is easily checked that $v$ is a corner of $U_{\tilde{T}}$ if and
only if $U_{\tilde{T}}$ contains unequal numbers of  black squares and white
squares adjacent to $v$ (and similarly for $U_{S(\tilde{T})}$). A domino $d \in
\tilde{T}$ may cover, of the four squares adjacent to $v$, one black square,
one white square, or one square of each color. In these three cases, $S(d)$
covers one white square, one black square, or no squares at all, respectively.
Thus the even vertex $v$  could be a corner of $U_{S(\tilde{T})}$ only if it
was already a corner of $U_{\tilde{T}}$. But we assumed $\tilde{T}$ was
odd-deficient, so that its only corners were at odd vertices. \qed

\vspace{0.2in}

Assume now that $\tilde{T}$ is an odd-deficient tiling, not of the entire
plane, but  of the order-$(n-1)$ Aztec diamond. We can use the above to show
that $S(\tilde{T})$ is an odd-deficient tiling of the order-$n$ diamond. It is
clear that for every domino $d \in \tilde{T}$, $S(d)$ lies in the order-$n$
diamond; what is less pictorially obvious is that the complement of
$S(\tilde{T})$ relative to the order-$n$ diamond must be a union of odd blocks.
One way to see this is to tile the complement of the order-$(n-1)$ Aztec
diamond in the fashion of Figure~\TiCo, obtaining an odd-deficient tiling
$\tilde{T}^+$ of the entire plane. Then by the Lemma, $S(\tilde{T}^+)$ is an
odd-deficient tiling of the plane; some of its missing odd-blocks lie in two
semi-infinite strips of height 2 to the left and right of the order-$n$
diamond, and all the others must lie strictly inside the order-$n$ diamond.
None of these blocks cross the boundary of the order-$n$ diamond, so if we add
these blocks to $S(\tilde{T})$, we get a complete tiling of the order-$n$
diamond.

Consider now an odd-deficient tiling $\tilde{T}$ of the order-$(n-1)$ Aztec
diamond, with $\tTvert$ equal to the set of vertical tiles of $\tilde{T}$; let
$$v(\tilde{T}) = 
\sum_{d \in \tilde{T}} v(d) =
\sum_{d \in \tTvert} v(d)$$
and
$$r(\tilde{T}) = 
\sum_{d \in \tilde{T}} r(d) =
\sum_{d \in \tTvert} r(d),$$
(recall that $v(d) = r(d) = 0$ for all horizontal dominoes $d$). Let
$$\AD (n-1,\tilde{T};x,q) = \sum x^{v(T)} q^{r(T)},$$
where the sum is over all tilings $T$ that extend $\tilde{T}$; we have
$$\AD (n-1;x,q) = \sum_{\tilde{T}} \AD (n-1,\tilde{T};x,q),$$
where the sum is over all partial tilings $\tilde{T}$ of the order-$(n-1)$
Aztec diamond. Say that $\tilde{T}$ is missing $m$ odd blocks, so that it gives
rise to $2^m$ distinct complete tilings $T$; then it is easily seen that
\begin{equation}
\label{easily}
\AD (n-1,\tilde{T};x,q) = 
(1+xq^{-1})^m \prod_{d \in \tilde{T}} x^{v(d)} q^{r(d)} \ .
\end{equation}
$S(\tilde{T})$ is an odd-deficient tiling of the order-$n$ Aztec diamond with
its odd coloring, missing $m+n$ odd blocks. Therefore, relative to the even
coloring, we have
$$\AD (n,S(\tilde{T});x,q)
= (1+xq)^{m+n} \prod_{d \in \tTvert} x^{v(S(d))} q^{-r(S(d))}.$$
The product in the right hand side can be rewritten as
\begin{eqnarray*}
\prod_{d \in \tTvert} x^{v(S(d))} q^{-r(S(d))}
& = & \prod_{d \in \tTvert}   x^{v(d)} q^{r(d)+1} \\
& = & \prod_{d \in \tTvert}   (xq^2)^{v(d)} q^{r(d)} \\
& = & \prod_{d \in \tilde{T}} (xq^2)^{v(d)} q^{r(d)} .
\end{eqnarray*}
But substituting $n$ for $n-1$ and $xq^2$ for $x$ in (\ref{easily}) yields
$$\AD (n,\tilde{T};xq^2,q) = (1+xq)^m 
\prod_{d \in \tilde{T}} (xq^2)^{v(d)} q^{r(d)} .$$
Hence
$$\AD (n,S(\tilde{T});x,q) = (1+xq)^n \AD (n-1,\tilde{T};xq^2,q).$$
Since every odd-deficient tiling of the order-$n$ Aztec diamond with odd
coloring is of the form $S(\tilde{T})$ for some odd-deficient tiling of the
order-$(n-1)$ Aztec diamond with even coloring, we can sum both sides of the
preceding equation over all $\tilde{T}$, obtaining
$$\AD(n;x,q) = (1+xq)^n \AD(n-1;xq^2,q).$$
The general formula for $\AD(n;x,q)$ follows immediately by induction.

Although this proof made no mention of alternating sign matrices, they are very
much involved in determining the exact locations of the various 2-by-2 blocks.
Specifically, let $T$ be a domino tiling of the Aztec diamond of order $n-1$,
and let $A$ be the $(n-1)$-by-$(n-1)$ alternating sign matrix determined by $T$
as in section 3. Then the locations of the odd blocks in $\tilde{T}$ are given
by the $1$'s in $A$, while the locations of the odd blocks in $S(\tilde{T})$
are given by the $-1$'s.

Latent within the proof of the formula for $\AD(n;x,q)$ is an iterative
bijection between domino-tilings of the order-$n$ Aztec diamond and bit-strings
of length $n(n+1)/2$. Say we are given a bit-string of length $1 + 2 + ... +
n$, and suppose we have already used the first $1 + 2 + ... + (k-1)$ bits to
construct a domino-tiling of the order-$(k-1)$ diamond. Impose the even
coloring on this Aztec diamond and locate the odd blocks, of which there are
$m$. Pick up these odd blocks in some definite order (of which we will say more
shortly) and put them elsewhere, retaining their order. Shuffle the dominoes in
the remaining partial tiling of the Aztec diamond of order $k-1$. The resulting
partial tiling of the order-$k$ Aztec diamond has $m+k$ holes in it; fill these
holes (again in some definite order) with the $m$ blocks that were removed
before, followed by $k$ other blocks, whose orientations (horizontal vs.\
vertical) are determined by the next $k$ bits of the bit-string. In this way
one obtains a complete tiling of the Aztec diamond of order $k$. Note that no
information has been lost; the procedure is fully reversible. Thus, iteration
of the process gives a bijection between bit-strings of length $n(n+1)/2$ and
domino tilings of the order-$n$ Aztec diamond. Moreover, every 0 (resp.\ 1) in
the bit-string leads to the creation of two horizontal (resp.\ vertical)
dominoes in the tiling, so it is immediate that the number of tilings of the
Aztec diamond with $2v$ vertical dominoes is ${n(n+1)/2} \choose v$.

The preceding construction requires a pairing between the $m$ missing odd
blocks of an odd-deficient tiling of the order-$(k-1)$ Aztec diamond and $m$ of
the $m+k$ missing odd blocks of an odd-deficient tiling of the order-$n$ Aztec
diamond. There is a canonical way of doing this pairing. Recall that these two
kinds of blocks correspond to the $-1$'s and $+1$'s in an alternating sign
matrix $A$, so it suffices to decree some sort of pairing between the $-1$'s
and a subset of the $+1$'s (which will leave $n$ $+1$'s left over). But this is
easy: just pair each $-1$ with the next $+1$ below it in its column. In terms
of shuffling, this means that the odd blocks of $\tilde{T}$ drift southeast
until they find a hole in $S(\tilde{T})$ that they can fit; this leaves $n$
holes near the upper left border of the order-$n$ Aztec diamond, which the $n$
new 2-by-2 blocks exactly fill.

It would be nice to have a ``shuffling'' proof of the general formula
(\ref{genform}) proved in section 4, and/or a procedure for randomly generating
monotone triangles according to the (uneven) probability distribution given by
the weights $W(\cdot)$.

\section{Square ice}

It is worthwhile to point out a connection between the combinatorial objects
investigated in this paper and a statistical mechanical model that has been
studied extensively since the 1960's. Recall that an $n$-by-$n$ alternating
sign matrix can be represented by its skewed summation, as in Figure~\Ice(a).
Replace each entry in the matrix by a node, and put a directed edge between
every two adjacent entries, pointing from the smaller to the larger. Then one
has a directed graph in which the circulation around every square cell is 0
(that is, each cell has two clockwise edges and two counterclockwise edges);
see Figure~\Ice(b). Finally, rotate each of these edges $90^\circ$
counterclockwise about its midpoint. The end result is a configuration like the
one shown in Figure~\Ice(c), with divergence 0 at each node (that is, each node
has two incoming arrows and two outgoing arrows). This is exactly the
square-ice model of statistical mechanics, with the special boundary condition
of incoming arrows along the left and right sides, and outgoing arrows along
the top and bottom. (For discussion of this and related models, see [\Ba] and
[\Pe].)

In the general square ice model, one associates a Boltzmann weight $\omega_i$
($i=1$ to 6) with each of the six possible vertex-configurations shown in
Figure~\VeCo; then the weight of a configuration is defined as
${\omega_1}^{k_1} {\omega_2}^{k_2} ...  {\omega_6}^{k_6}$, where $k_i$ is the
number of vertices in the lattice of type $i$, and the {\bf partition function}
associated with the model (denoted by $Z$) is the sum of the weights of all
possible configurations. $Z$ has an implicit dependence on the lattice-size and
the boundary conditions. It is customary to impose periodic boundary
conditions, but we instead impose the ``in-at-the-sides,
out-at-the-top-and-bottom'' condition on our $n$-by-$n$ grid. Call this the
{\bf Aztec boundary condition}.

To recast our work on domino tilings of the Aztec diamond in terms of square
ice, it is convenient to rephrase domino-tilings as dimer arrangements, or
1-factors. Specifically, we define a graph $G'$ whose vertices correspond to
the cells of the order-$n$ Aztec diamond, with an edge between two vertices of
$G\,'$ if and only if the corresponding cells are adjacent. Then a
domino-tiling of the Aztec diamond corresponds to a 1-factor $F$ of $G'$ (a
collection of disjoint edges covering all vertices).

There is a general method for writing the number of 1-factors of a planar
graph  as a Pfaffian ([\Kab]). Indeed, if one assigns weight $w(e)$ to each
edge $e$ of a planar graph on $N$ vertices and defines the weight of a 1-factor
as the product of its constituent weights, then the sum of the weights of all
1-factors of the graph is equal to the Pfaffian of an antisymmetric $N$-by-$N$
matrix whose $i,j$th entry is $\pm w(e)$ if the graph has an edge $e$ between
$i$ and $j$ and 0 otherwise. (The delicate point is the correct choice of
signs.) This method has been applied to the problem of counting 1-factors of
$m$-by-$n$ grids (equivalently, domino-tiling of $m$-by-$n$ rectangles); see
[\Kaa], [\Be], [\Lo]. The Pfaffian method provides a yet another route to our
result on tilings of the Aztec diamond, though we have omit the calculation
here; see [\Ya].

It is convenient to rotate the graph $G'$ $45^\circ$ clockwise, as in
Figure~\Tilt(a). Call a cell of $G'$ even or odd according to the parity of the
corresponding vertex of $G$ (under the standard coloring), so that the four
extreme cells of $G'$ are even. Every even cell is bounded by four edges, of
which two, one, or none may be present in any particular 1-factor; the seven
possibilities appear at the top of Figure~\Trans, where a bold marking
indicates the presence of an edge. If we replace each even cell by the
corresponding ice-junction given at the bottom of the Figure~\Trans, it is easy
to check that the result is a valid ice-configuration satisfying our special
boundary conditions, and that every such configuration arises in this way. The
process is exemplified in Figure~\Tilt(b). Note that the transformation from
1-factors to ice-configurations is not one-to-one; it is in fact
$2^{k_5}$-to-one, where $k_5$ is the number of vertices of type 5 in the ice
pattern. It can be checked that this transformation is equivalent to the more
roundabout operation of converting the 1-factor to a domino-tiling, using the
heights of the even vertices to form an $n$-by-$n$ alternating sign matrix, and
then turning the matrix into an ice-pattern as in the first paragraph of this
section.

Let $T$ be a tiling of the Aztec diamond, and let $F$ be the associated
1-factor of $G'$. Note that every domino in $T$ corresponds to an edge in $F$,
and that this edge belongs to a unique even cell of $G'$. Hence, if we assign
the weights $x$, $x$, $1$, $1$, $1$, $x^2$, and $1$ to the respective
cell-figures, the product of the weights of the cell-figures appearing in $F$
is equal to $x^{2v(T)}$. Thus, if we put
\begin{eqnarray*}
\omega_1 & = & x, \\
\omega_2 & = & x, \\
\omega_3 & = & 1, \\
\omega_4 & = & 1, \\
\omega_5 & = & 1+x^2, \ \mbox{and} \\
\omega_6 & = & 1,
\end{eqnarray*}
then the partition function $Z$ coincides with the generating function
$$\AD(n;x^2) = (1+x^2)^{n(n+1)/2}.$$

Note that $k_5 - k_6 = n$  for all order-$n$ ice-configurations with Aztec
boundary condition (corresponding to the fact that the number of $1$'s in an
$n$-by-$n$ alternating sign matrix must be $n$ more than the number of $-1$'s).
Hence replacing $\omega_5$ and $\omega_6$ by $\sqrt{1+x^2}$ merely divides the
partition function by $(1+x^2)^{n/2}$. Furthermore, $k_1 + k_2 + ... + k_6 =
n^2$, so multiplying all the Boltzmann weights by a factor $b$ merely
multiplies the partition function by $b^{n^2}$. Writing $a = bx$ and $c =
b\sqrt{x^2+1} = \sqrt{a^2+b^2}$, we see (after an easy calculation) that for
the square-ice model  with Aztec boundary condition and with Boltzmann weights
\begin{eqnarray*}
\omega_1 & = & a, \\
\omega_2 & = & a, \\
\omega_3 & = & b, \\
\omega_4 & = & b, \\
\omega_5 & = & c, \\
\omega_6 & = & c
\end{eqnarray*}
satisfying $a^2 + b^2 = c^2$, the partition function is given by $Z = c^{n^2}$.

It should be noticed that this family of special cases of the ice model (given
by $a,b,c$ satisfying $a^2+b^2=c^2$) is also the family that corresponds to the
free fermion case, and is precisely the case in which the model has been solved
by the method of Pfaffians. (See p.\ 151, 270-271 of [\Ba], as well as [\FW].)
This leads us to suspect that domino shuffling may in fact arise from some
combinatorial interpretation of the Pfaffian solution.

We must emphasize the role played by the Aztec boundary conditions in the
foregoing analysis, since it adds an element essentially foreign to the
physical significance of the ice model. In particular, Lieb's solution of the
ice model in the case $\omega_1 = \omega_2 = ... = \omega_6 = 1$ [\Li] tells us
that there are asymptotically
$$\sqrt{64/27} ^ {\ n^2}$$
order-$n$ ice-configurations with periodic boundary conditions; on the other
hand, if the conjecture of Mills, Robbins, and Rumsey is correct, the number of
order-$n$ ice-configurations with Aztec boundary conditions should
asymptotically be only
$$\sqrt{27/16} ^ {\ n^2}.$$
Clearly there are more constraints on a domino tiling near the boundary of an
Aztec diamond than there are near the middle; this accounts for at least some
of the drop in entropy. It would be interesting to know in a more quantitative
way how the entropy of a random tiling is spatially distributed throughout a
large Aztec diamond.

\section{Epilogue}

There have been many combinatorial transformations in this article,
so it may be useful to review them.

First, we have:
\begin{enumerate}
\item[(1)] tilings; 
\item[(2)] height-functions associated with tilings; and
\item[(3)] the order ideals associated with those height-functions.
\end{enumerate}
We saw how to go from 1 to 2 (Thurston's marking scheme),
from 2 to 3 (see the construction of the poset $P$ in Section 3),
and from 3 back to 1 (the stacked cubes).

Then we have:
\begin{enumerate}
\item[(4)] alternating sign matrices;
\item[(5)] height-functions associated with alternating sign matrices;
\item[(6)] the order ideals associated with those height-functions;
\item[(7)] monotone triangles;
\item[(8)] states of the square ice model (or equivalently its dual).
\end{enumerate}
We saw the correspondence between (4) and (5) and between (5) and (6) in
Section 3, between (4) and (7) in Section 4, and between (5) and (8) in Section
7. Further correspondences can be made. For instance, to get from (4) to (8)
directly, we replace a 1 in an alternating sign matrix by a
vertex-configuration of type 5, a -1 by a by a vertex-configuration of type 6,
and each 0 by the unique vertex-configuartion of type 1-4 which fits in the
pattern (note that arrows ``go straight through'' configurations of type 1-4
without reversing).

Then there are the mappings between (1)-(3) and (4)-(8), under the
correspondence between domino tilings and compatible pairs of alternating sign
matrices. We saw in Section 3 how to pass between (2) and (5), and between (3)
and (6). Other connections can be made, and the reader might find it
instructive to try to establish them.

There are actually even more incarnations of alternating sign matrices than
have been discussed here: 3-colorings of certain graphs (subject to boundary
constraints), 2-factors of some related graphs, and tilings of various regions
in the plane by shapes of two kinds. These other structures may be discussed in
a future paper. Then there are other combinatorial objects which appear (but
have not been proved) to be equinumerous with the alternating sign matrices,
namely, descending plane partitions and self-complementary totally symmetric
plane partitions.  See [\Ro] for details.

Richard Stanley has discovered that our two-variable generating function for
tilings of the order-$n$ Aztec diamond is actually a specialization of a
$2n$-variable generating function. A proof of this identity via the shuffling
method of Section 6 is described in [\Ya].

\section*{Acknowledgments}

We wish to thank Carol Sandstrom and Chris Small for suggesting some of the
terminology used in this article.
 
\section{Bibliography}

\begin{description}

\item[{\bf [\Ba]}]
R.\ J.\ Baxter, {\it Exactly Solved Models in Statistical Mechanics.}
Academic Press, 1982.

\item[{\bf [\Be]}]
E.\ F.\ Beckenbach, ed., {\it Applied Combinatorial Mathematics.}
J.\ Wiley, 1964.

\item[{\bf [\CL]}]
J.\ H.\ Conway and J.\ C.\ Lagarias,
Tilings with polyominoes and combinatorial group theory,
{\it J.\ Combin.\ Theory A} {\bf 53} (1990), 183-208.

\item[{\bf [\FW]}]
C.\ Fan and F.\ Y.\ Wu, General lattice model of phase transitions,
{\it Phys. Rev.} {\bf B2} (1970), 723-733.

\item[{\bf [\Gr]}]
J.\ A.\ Green, {\it Polynomial Representations of $GL_n$.}
Springer Lecture Notes \#830, 1980.

\item[{\bf [\Kaa]}]
P.\ W.\ Kasteleyn, The statistics of dimers on a lattice,
I: The number of dimer arrangements on a quadratic lattice,
{\it Physica} {\bf 27} (1961), 1209-1225.

\item[{\bf [\Kab]}]
P.\ W.\ Kasteleyn, 
Graph theory and crystal physics, in
F.\ Harary, ed.,
{\it Graph Theory and Theoretical Physics.}
Academic Press, 1967.

\item[{\bf [\Li]}]
E.\ Lieb, Residual entropy of square ice,
{\it Phys.\ Rev.} {\bf 162} (1967), 162-172.

\item[{\bf [\Lo]}]
L.\ Lov\'asz, {\it Combinatorial Problems and Exercises,}
problem 4.29.
North Holland, 1979.

\item[{\bf [\MRRa]}]
W.\ H.\ Mills, D.\ P.\ Robbins, and H.\ Rumsey, Jr.,
Proof of the Macdonald conjecture,
{\it Invent.\ Math.} {\bf 66} (1982), 73-87.

\item[{\bf [\MRRb]}]
W.\ H.\ Mills, D.\ P.\ Robbins, and H.\ Rumsey, Jr.,
Alternating sign matrices and descending plane partitions,
{\it J.\ Combin.\ Theory A} {\bf 34} (1983), 340-359.

\item[{\bf [\Pe]}]
J.\ K.\ Percus,
{\it Combinatorial Methods.}
Courant Institute, 1969.

\item[{\bf [\PS]}]
G.\ P\'olya and S.\ Szeg\"o,
{\it Problems and Theorems in Analysis, Volume II.}
Springer, 1976, p.\ 134.

\item[{\bf [\Ro]}]
D.\ P.\ Robbins,
The story of 1, 2, 7, 42, 429, 7436, ...
(preprint);
to appear in {\it The Mathematical Intelligencer.}

\item[{\bf [\RR]}]
D.\ P.\ Robbins and H.\ Rumsey, Jr.,
Determinants and alternating sign matrices,
{\it Adv.\ in Math.} {\bf 62} (1986), 169-184.

\item[{\bf [\Sp]}]
A.\ E.\ Spencer,
Problem E 2637,
{\it American Mathematical Monthly}, {\bf 84} (1977), 134-135;
solution published in {\bf 85} (1978), 386-387.

\item[{\bf [\Sta]}]
R.\ Stanley,
{\it Enumerative Combinatorics, Volume I.}
Wadsworth and Brooks/Cole, 1986.

\item[{\bf [\St]}]
R.\ Stanley,
A baker's dozen of conjectures concerning plane partitions
(preprint).

\item[{\bf [\Th]}]
W.\ Thurston,
Conway's tiling groups
(preprint);
to appear in {\it The American Mathematical Monthly.}

\item[{\bf [\Ya]}]
B.-Y.\ Yang,
Three Enumeration Problems Concerning Aztec Diamonds.
Submitted as an M.I.T.\ doctoral thesis, May 1991.

\end{description}

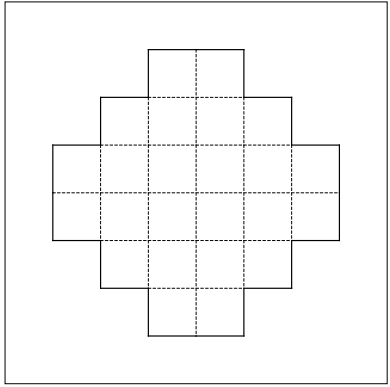
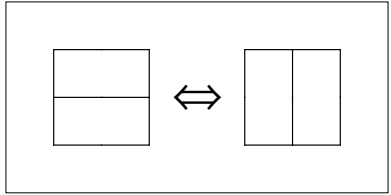
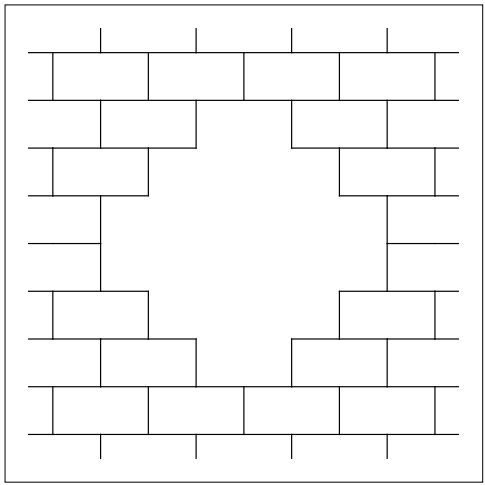
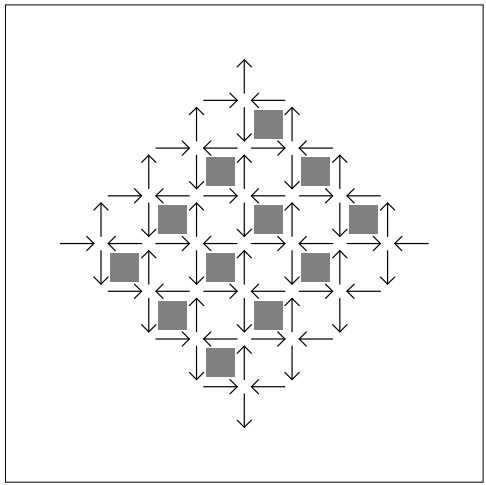
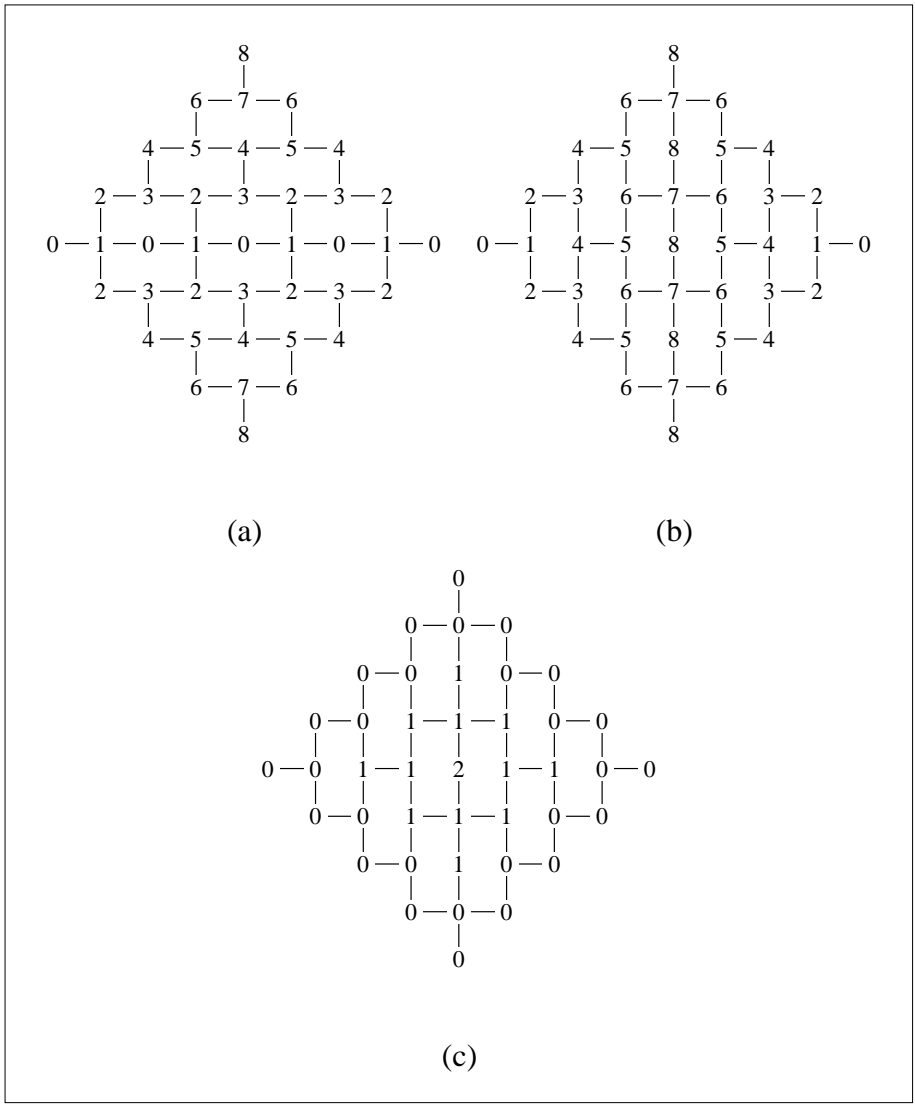
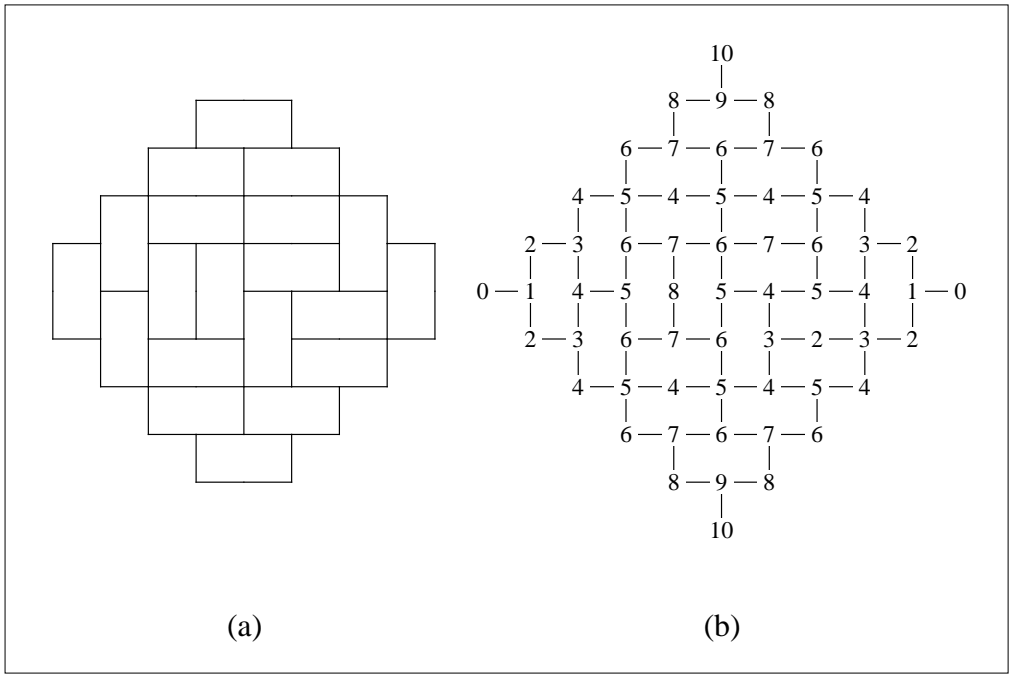
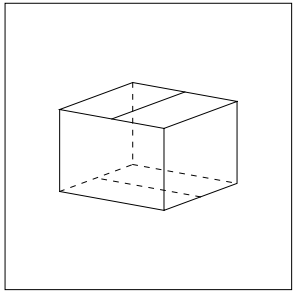
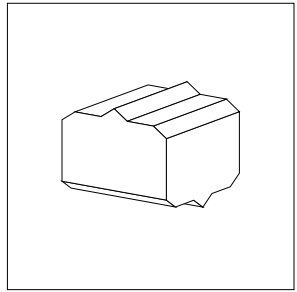
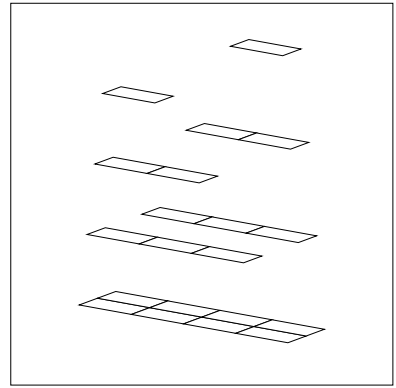
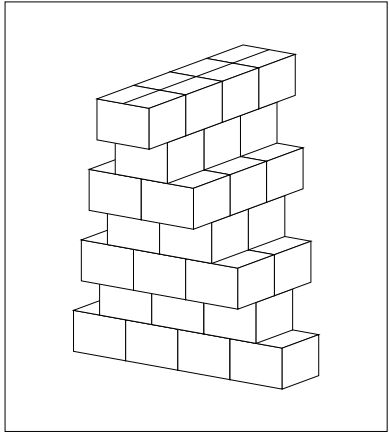
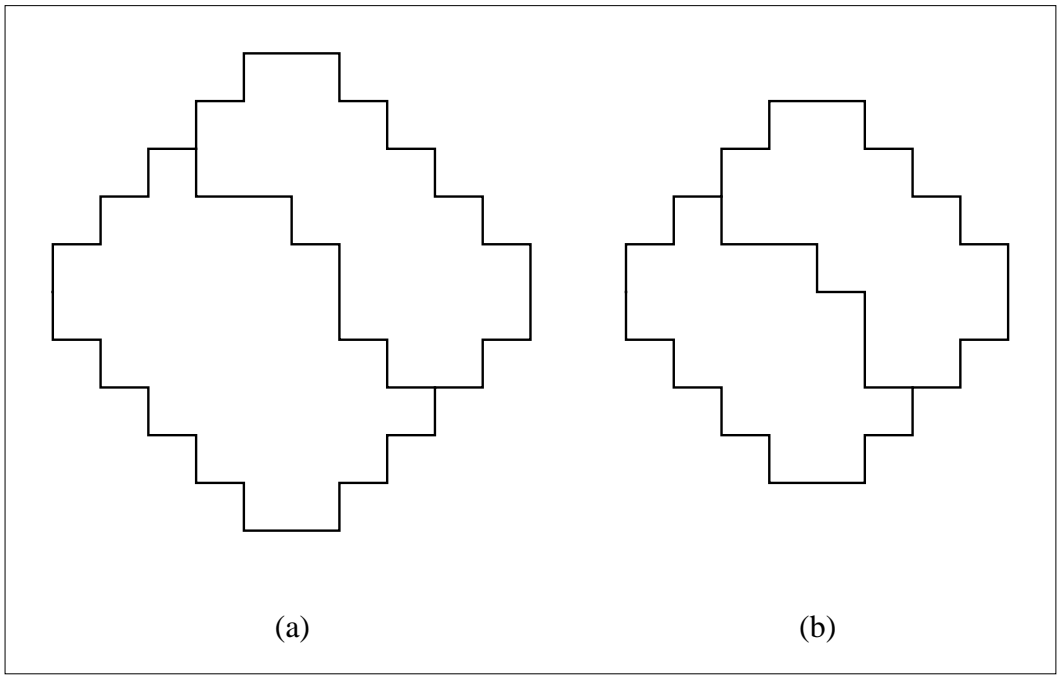
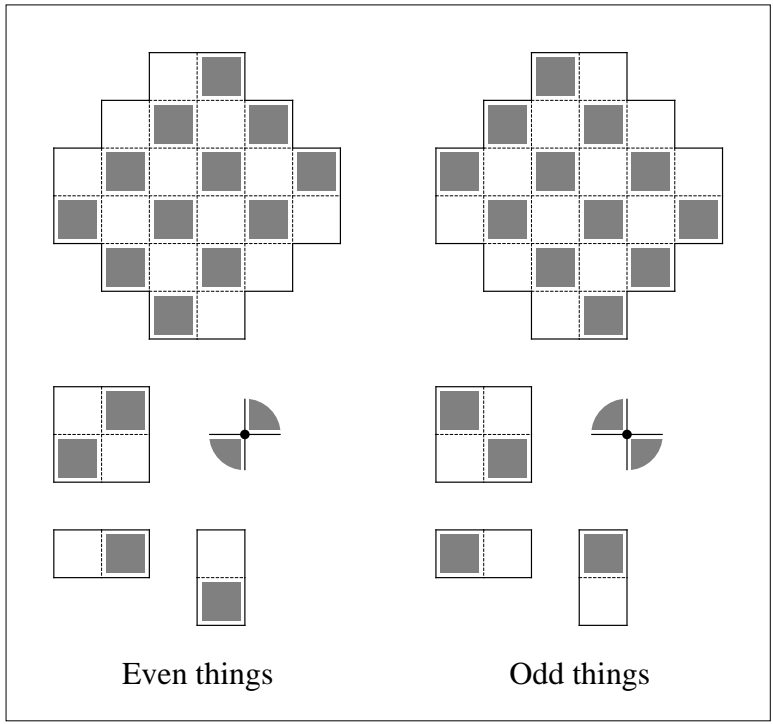
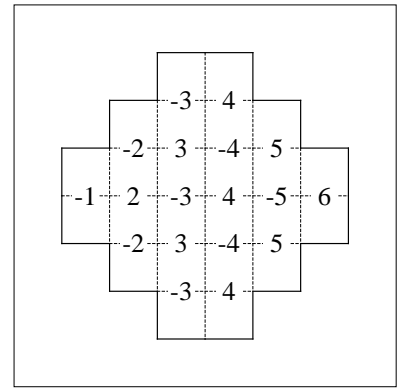
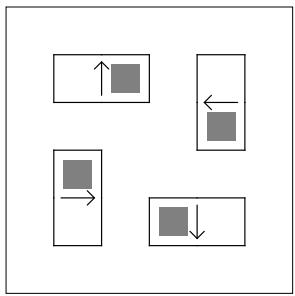
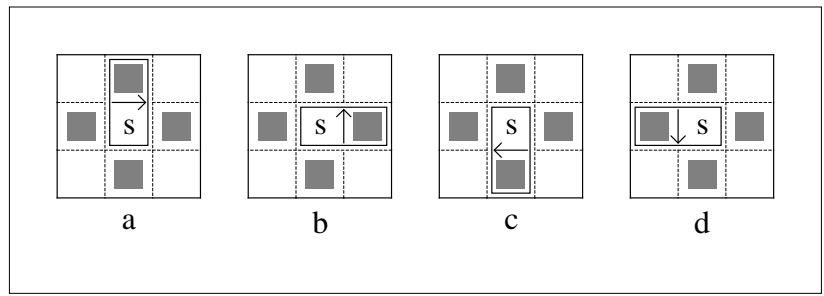
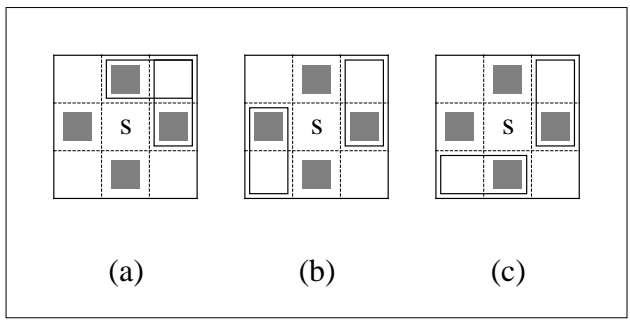
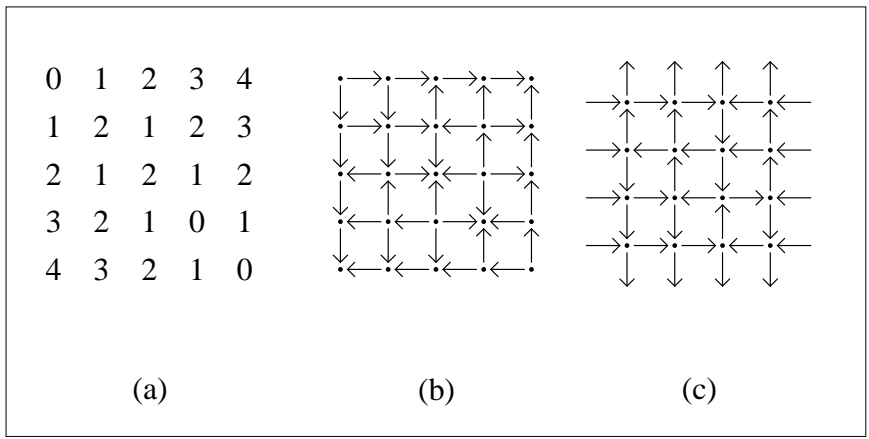
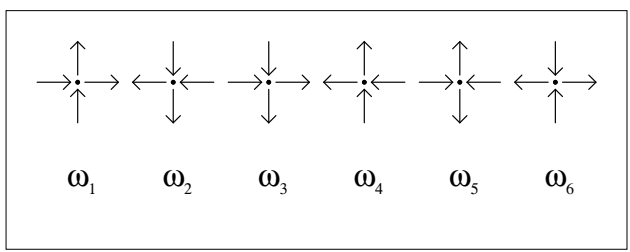
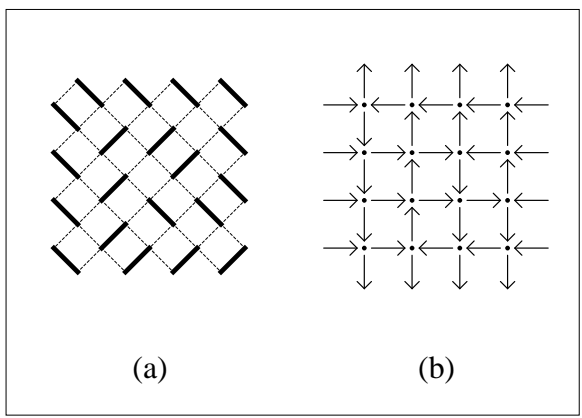
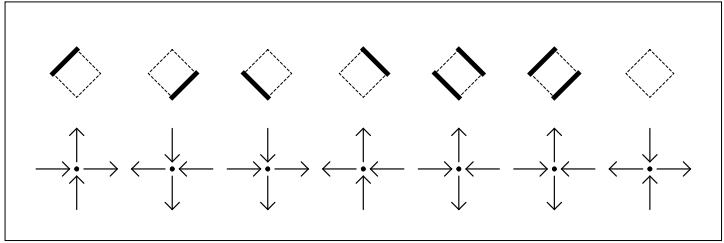
\end{document}